\documentclass[11pt]{article}
\usepackage{mathrsfs}

\usepackage{latexsym}
\usepackage{bbm}
\usepackage{amsmath}
\usepackage{amsfonts}
\usepackage{amssymb}
\usepackage{multirow}
\usepackage{latexsym}
\usepackage[dvips]{graphicx}
\usepackage{overpic}
\usepackage{graphicx}
\newtheorem{theorem}{\bf Theorem}[section]
\newtheorem{lemma}[theorem]{\bf Lemma}

\newtheorem{defn}{\bf Definition}[section]

\newenvironment{proof}{\noindent{\em Proof.}}{\quad \hfill$\Box$\vspace{2ex}}

\def \bZ {\Bbb Z}

\def \bR {\Bbb R}

\def \and {\, \mbox{\rm and}\, }

\def \supp {\,{\rm supp}\,}

\def \Re {\,{\rm Re}\,}
\def \Im {\,{\rm Im}\,}

\def \l {\left}
\def \r {\right}
\makeatletter

\newcommand{\Rmnum}[1]{\expandafter\@slowromancap\romannumeral #1@}
\makeatother

\begin{document}
\title{\bf Sharp $L^{p}$-Boundedness of Oscillatory Integral Operators with
 Polynomial Phases}
\author{Zuoshunhua Shi \thanks{School of Mathematical Sciences, University
 of Chinese Academy of Sciences, Beijing 100049, P. R. China. E-mail address:
{\it shizuoshunhua11b@mails.ucas.ac.cn}.} \quad and \quad Dunyan
Yan\thanks{School of Mathematical Sciences, University of Chinese Academy of Sciences, Beijing 100190, P. R. China. E-mail
address: {\it ydunyan@ucas.ac.cn}.} }
\date{}
\maketitle{}
\begin{abstract}
In this paper, we shall prove the $L^{p}$ endpoint decay estimates of oscillatory integral operators with homogeneous polynomial phases $S$ in $\mathbb{R} \times \mathbb{R}$. As a consequence, sharp $L^{p}$ decay estimates are also obtained when polynomial phases have the form $S(x^{m_{1}},y^{m_{2}})$ with $m_1$ and $m_2$ being positive integers .
\end{abstract}
\section{Introduction}

\hspace{.5cm}In this paper, we mainly consider the following operator
\begin{equation}\label{Osc integral operator with cut-off}
T_{\lambda}(f)(x)=\int_{-\infty}^{\infty}e^{i\lambda S(x,y)}\varphi
(x,y)f(y)dy,
\end{equation}
where $\lambda \in \mathbb{R}$, $\varphi \in
C_{0}^{\infty}(\mathbb{R}^{2})$ and $S$ is a real-valued homogeneous
polynomial in $\mathbb{R}\times \mathbb{R}$. We can write the phase
$S$ as
\begin{equation}\label{homogeneous polynomial}
S(x,y)=\sum_{k=0}^{n}a_{k}x^{n-k}y^{k}
 \end{equation}
with real coefficients $a_{k}$.

If $S$ is a general real-valued smooth function with $S_{xy}''\neq 0$ on the support of the cut-off function, H\"{o}rmander (\cite{hormander2}) obtained the sharp $L^2$ operator norm estimate  $\|T_{\lambda}\|\leq C|\lambda|^{-1/2}$. If the phase is degenerate on $\supp(\varphi)$, then sharp decay estimates cannot be obtained directly and we need a suitable resolution of the singular variety $\{S_{xy}''=0\}$. In this direction, the scalar oscillatory integrals with analytic phases were studied by Varchenko \cite{varchenko}. For homogeneous polynomials $S$, the study of these operators was initiated from \cite{PS1992}. In \cite{PS1994}, Phong and Stein gave a necessary and sufficient condition under which the $L^2$ operator norm satisfies the sharp estimate $\|T_{\lambda}\|\leq C|\lambda|^{-\frac1n}$. For general real analytic phases, they established the relation between the sharp $L^2$ estimate and the Newton polyhedron of $S$; see \cite{PS1997}. In \cite{rychkov}, Rychkov extended $L^2$ estimates in \cite{PS1997} to most smooth phases and full generalizations to smooth phases were proved by Greenblatt in \cite{greenblatt2}. For other related results, we refer the reader to \cite{greenblatt1}, \cite{greenleafseeger1}, \cite{seeger1} and \cite{GPT}.

Some averaging operators of Radon transforms are closed related to above oscillatory integral operators. Denote $R$ by the Radon transform
\begin{equation*}
Rf(x)=\int_{-\infty}^{\infty} f(x_1+S(x_2,t),t)\varphi(x,t)dt
\end{equation*}
with $\varphi\in C_0^\infty(\mathbb{R}^2\times\mathbb{R})$. The sharp $L^p-L^q$ estimates and $L^p$ Sobolev regularity were obtained by Phong and Stein in \cite{PS1994} for homogeneous polynomials $S$ except endpoint estimates. For analytic phases, endpoint $L^p-L^q$ estimates are previously known and sharp $L^p$ Sobolev regularities were obtained except extreme points; see \cite{lee} and \cite{yangchanwoo2} as well as \cite{bak} and \cite{bos} by imposing certain left and right finite type conditions. It is notable that endpoint $L^p$ Sobolev regularity may fail; see \cite{christ}. For more general Radon transforms which can regarded as degenerate Fourier integral operators, we refer the reader to the survey paper \cite{greenleafseeger2} and references therein.

Our objective is to establish the sharp $L^p$ estimates for $T_\lambda$. In the case of two sided fold singularities, Greenleaf and Seeger obtained in \cite{greenleafseeger1} the sharp $L^p$ estimates for a smooth phase $S$. The sharp $L^p$ estimates were obtained in \cite{yangchanwoo} for homogeneous polynomial phases $S$ under the assumption $a_1a_{n-1}\neq 0$. For an analytic phase $S$, sharp $L^p$ estimates of $T_\lambda$ had been established in \cite{yangchanwoo2} except some extreme points of the reduced Newton polyhedron of $S$. Assume $S$ is a real analytic function near the origin. Then the power expansion gives $S(x,y)=\sum_{k,l\geq0}a_{k,l}x^ky^l$ near the origin. The Newton polyhedron is the convex hull of the sets $\{(x,y):x\geq k, y\geq l\}$ with nonnegative integers $k$ and $l$ satisfying $a_{k,l}\neq 0$. Similarly, the reduced Newton polyhedron $\mathcal{N}(S)$ is obtained by taking the same convex hull with the additional condition $kl\neq 0$. Let $T_\lambda$ be defined as in (\ref{Osc integral operator with cut-off}) with the cut-off function $\varphi$ supported in a small neighborhood of the origin. If $k$ and $l$ are two positive integers such that $|\partial_x^k\partial_y^lS(0,0)|\neq 0$, then the endpoint $L^p$ estimates are given by
\begin{equation}\label{sharp estimate for analytic phase}
\|T_\lambda\|_{L^{(k+l)/k}\rightarrow L^{(k+l)/k}}\leq C|\lambda|^{-1/(k+l)}.
\end{equation}
By interpolation, it is easily verified that these estimates imply the sharp $L^2$ decay rate obtained by Phong and Stein in \cite{PS1997}. We shall point out that the above estimates are sharp provided that $(k,l)$ is a vertex of the reduced Newton polyhedron of $S$; see \cite{yangchanwoo2}. If $(k,l)$ is not a vertex but lying on the boundary of $\mathcal{N}(S)$, the above sharp estimates were obtained in \cite{yangchanwoo2}. But only weak type results were proved in \cite{yangchanwoo2} when $(k,l)$ is a vertex. The question arises naturally whether the inequality (\ref{sharp estimate for analytic phase}) is true for extreme points $(k,l)$. In this paper, we shall give an affirmative answer when $S$ is a real-valued homogeneous polynomial. Sharp estimates (\ref{sharp estimate for analytic phase}) are also true for a more general class of polynomials.

Now we turn to the issue for homogeneous phases $S$ as in (\ref{homogeneous polynomial}). To avoid triviality, we assume $a_k\neq 0$ for some
$1\leq k \leq n-1$. Define $k_{min}$ and $k_{max}$ as follows
\begin{equation}\label{k_min}
k_{min}=\min\big\{1 \leq k \leq n-1: a_{k}\neq 0\big\}
\end{equation}
and
\begin{equation}\label{k_max}
 k_{max}=\max\big\{1 \leq k \leq n-1: a_{k}\neq 0\big\}.
\end{equation}
It is clear that $(k_{min},n-k_{min})$ and $(k_{max},n-k_{max})$ are two vertices of $\mathcal{N}(S)$.

Combining previously known results, we can state our main result as follows. We emphasize that only the endpoint estimates are new.

\begin{theorem}\label{main theorem 1}
Suppose that $T_\lambda$, $S$, $k_{min}$ and $k_{max}$ are given as
above. Then the decay estimate
\begin{equation}\label{main decay estimate}
\|T_\lambda\|_{L^p\rightarrow L^p} \leq C|\lambda|^{-\frac1n}
\end{equation}
holds for any amplitude $\varphi \in C_{0}^\infty$ if and only if $\frac{n}{n-k_{min}}\le p \le \frac{n}{n-k_{max}}$.
\end{theorem}

At the same time, we also consider the following two operators
\begin{equation}\label{OIO with homogeneous polynomial without cut-off}
T(f)(x)=\int_{-\infty}^{\infty}e^{iS(x,y)}f(y)dy
\end{equation}
and
\begin{equation}\label{OIO with homogeneous polynomial without cut-off-1}
T_{m_{1},m_{2}}(f)(x)=\int_{-\infty}^{\infty}e^{iS(x^{m_1},y^{m_2})}f(y)dy
\end{equation}
for $f \in C_{0}^{\infty}$ and positive integers $m_1$ and $m_2$. By a routine scaling argument,  $L^p$-boundedness of $T$ is equivalent to the decay estimate (\ref{main decay estimate}). In the following sections, we shall prove Theorem \ref{main theorem 1} by showing that $T$ is bounded on $L^{p}$ for $p$ in the range described as above. As a consequence, we also obtain the sharp $L^p$
boundedness of $T_{m_1,m_2}$ by invoking a simple interpolation lemma.

\begin{theorem}\label{Lp result for OIO with nonhomogeneous phases}
Let $m_1$ and $m_2$ be two positive integers. Then $T_{m_{1},m_{2}}$ defined by {\rm(\ref{OIO with homogeneous polynomial without cut-off-1})}
 has a bounded extension from $L^{p}$ to itself if and only if
\begin{equation*}
\frac{k_{min}m_2}{(n-k_{min})m_1}+1\leq p \leq
\frac{k_{max}m_2}{(n-k_{max})m_1}+1.
\end{equation*}
\end{theorem}

Now we first present some previously known results. When $S(x,y)=c(x-y)^n$ for nonzero $c$,  the $L^p$ boundedness of $T$  was established by \cite{jurkatsampson} and \cite{sjolin}
when $(x-y)^n$ is replaced by $|x-y|^n$; see also \cite{lu-zhang}. The arguments in the previous papers are also applicable for $S(x,y)=c(x-y)^n$. Another simple case is $S(x,y)=a_{k}x^{n-k}y^{k}$ with $a_k\neq 0$ and the corresponding result is contained in \cite{pansampson}. If $S$ is not of the form $c(x-\alpha y)^n$, it was proved that (\ref{main decay estimate}) holds when $a_{1}a_{n-1}\neq 0$ in \cite{yangchanwoo}. For oscillatory integral operators with two sided fold singularities, we refer the reader to \cite{greenleafseeger1} for sharp $L^p$ estimates. In \cite{yangchanwoo2}, Yang obtained the sharp $L^p$ decay estimates except the extreme points of $\mathcal{N}(S)$ when $S$ is a real analytic phases.

Our proof of the main result relies on a complex interpolation between $H_E^1-L^1$ and $L^2-L^2$. This method appeared earlier in \cite{PS1986} and \cite{greenleafseeger1}. Here $H_E^1$ ia a variant of Hardy spaces which will be defined later. To obtain sharp $H_E^1-L^1$ and $L^2-L^2$ estimates, we shall exploit a family of damped oscillatory integral operators. By insertion of the damping factor $|S_{xy}^{''}|^{1/2}$ in (\ref{OIO with homogeneous polynomial without cut-off}),
a useful sharp estimate has been obtained in \cite{PS1994}; see also \cite{PS1998} for the treatment
of analytic phases. For convenience, we present the result in \cite{PS1994}.

\begin{theorem}\label{Lp boundedness with damping factor}{\rm(\cite{PS1994})}
Assume that $S(x,y)$ is a real-valued homogeneous polynomial as in
\textrm{\rm(\ref{homogeneous polynomial})} with $a_{k}\neq 0$ for
some $1\leq k \leq n-1$. Let $U$ be given by
$$U(f)(x)=\int_{-\infty}^{\infty}e^{iS(x,y)}|S''_{xy}(x,y)|^{\frac12}f(y)dy$$
for $f\in C_0^{\infty}$. Then $U$ extends as a bounded operator from
$L^2$ to itself.
\end{theorem}

It was remarked in \cite{PS1994} that the same result still  holds
if the damping factor $|S''_{xy}(x,y)|^{1/2}$ is replaced by
$|S''_{xy}(x,y)|^{\alpha}$ with $\Re(\alpha)=1/2$ for which the
operator norm is bounded by a constant multiple of $(1+|\Im(\alpha)|)^{2}$. This provides the sharp endpoint $L^2$ estimates. On the other hand, we shall see that the sharp $H_E^1-L^1$ estimate (without decay on the parameter $\lambda$) is closely related to a class of oscillatory singular integral operators considered by Ricci and Stein in \cite{ricci},
\begin{equation}\label{osc singular integral with polynomial phases}
T_{P}(f)(x)={\rm p.v.}\int_{\mathbb{R}^n}e^{i P(x,y)}K(x,y)f(y)dy,
\end{equation}
where $f$ is initially assumed to be smooth with compact support, $P(x,y)$ a real polynomial in $\mathbb{R}^n\times \mathbb{R}^n$ and $K(x,y)$ a Calder\'{o}n-Zygmund kernel. It is known that
$T_P$ extends as a bounded operator from $L^p$ to itself for $1<p<\infty$ and the operator norm depends only on the degree of $P$ but not on its coefficients; see \cite{ricci} and \cite{lu-zhang}. In the case $p=1$, $T_{P}$ is bounded from $H^1_{E}$ to $L^1$; see \cite{pan}. Now we define the space $H^1_E$, associated to a polynomial $P(x,y)$, using $H_{E}^1$ atoms (\cite{pan}).

\begin{defn}\label{definition of variants of H^1}
Suppose that $P(\cdot,\cdot)$ is a real-valued polynomial in
$\mathbb{R}^{n}\times \mathbb{R}^{n}$. Let $Q\subset \mathbb{R}^{n}$
be a cube with sides parallel to the axes and the center $x_{Q}$.
An $H^{1}_{E}$ atom associated to the polynomial $P$ is a
measurable function $a$ satisfying\\

{\rm{(i)}} {\rm \supp($a$)} $\subset$ $Q${\rm;}~~~~~{\rm{(ii)}} $|a(x)| \leq |Q|^{-1}$, a.e. $x\in Q${\rm ;}~~~~~{\rm{(iii)}} ${\displaystyle \int e^{iP(x_{Q},y)}a(y)dy=0}$.
\end{defn}

Then $H^1_E(\mathbb{R}^n; P)$ consists of all those
$f\in L^{1}(\mathbb{R}^n)$ which can be decomposed as
$f=\sum_{j}\lambda_{j}a_{j},$
where each $a_{j}$ is an atom associated to $P$ and
$\lambda_{j}\in \mathbb{C}$ satisfying $\sum_{j}|\lambda_{j}|<\infty$.
The norm of $f$ in $H^{1}_{E}$ is given by
$\|f\|_{H^{1}_{E}}
=\inf\l\{\sum|\lambda_j|: \;f=\sum \lambda_j a_j,
\; a_{j}\in H_{E}^1\r\}.$

In the theory of singular integral operators, it is well known that the spaces $H^{1}$ and $BMO$ are appropriate substitutes of $L^{1}$ and $L^{\infty}$, respectively. One aspect of this is that they
play an important role in the interpolation of operators. A useful device is the sharp function invented by C. Fefferman and E. M. Stein in \cite{feffermanstein}. The dual space $BMO_E$ of $H^{1}_E$ and the associated sharp function are defined as follows (see also \cite{pan} and \cite{PS1986}).

\begin{defn}\label{definition of variants of BMO}
Let $f^{\sharp}_{E}$ be the sharp function given by
\begin{equation}\label{sharp function }
f^{\sharp}_{E}(x)=\sup_{Q \ni x}\frac{1}{|Q|}\int_Q
\l|f(y)-f_{Q}^{E}(y)\r|dy
\end{equation}
with
$$f_{Q}^{E}(x)=e^{iP(x_Q,x)}\frac{1}{|Q|}\int_Q e^{-iP(x_Q,y)}f(y)dy.$$
Then $BMO_{E}$, associated to the polynomial $P$, consists of all
locally integrable functions $f$ such that
$\|f\|_{BMO_E}=\|f^{\sharp}_E\|_{L^{\infty}}<\infty$.
This finite number is defined to be the norm of $f$ in $BMO_E$.
\end{defn}

The present paper is organized as follows. $\S$2 contains some basic
lemmas and section 3 is devoted to mapping properties of some
fractional integral operators. In $\S$4, we shall treat oscillatory
integral operators with polynomial phases and prove the estimate
$H_E^1-L^1$. The sharp $L^2$ decay estimates are obtained in $\S$5
for damped oscillatory integral operators. The main result of
Theorem \ref{main theorem 1} is proved in $\S$6. In the final part
$\S$7, we shall apply Theorem \ref{main theorem 1} to prove Theorem
\ref{Lp result for OIO with nonhomogeneous phases} and give an easy
proof of Pitt's inequality. The symbol $C$ stands for a constant
which may vary from line to line.

\section{Some Basic Lemmas}
The following lemma is useful in our proof of Theorem \ref{Lp result for OIO with nonhomogeneous phases}.
\begin{lemma}\label{interpolation lemma}
Let $S$ be a sublinear operator which is initially defined for
simple functions in $\mathbb{R}$. If there exists a constant $C>0$
such that~
\textrm{\rm{(i)}} $\|Sf\|_{\infty}\leq C\|f\|_{1}$ and
\textrm{\rm{(ii)}} $\|Sf\|_{p_{0}}\leq C \|f\|_{p_{0}}$
for some $1<p_0<\infty$ and all simple functions $f$, then the following inequality
$$\int_{\bR}|Sf(x)|^{p}|x|^{(p-p_{0})/(p_{0}-1)}dx \leq C \int_{\bR} |f(x)|^{p}dx$$
holds for $ 1< p \leq p_0$, where  the constant $C$ is independent of
$f$.
\end{lemma}
For the special case $p_{0}=2$, this lemma is contained in
\cite{pansampson}.

\begin{proof}
Let
$Tf(x)=|x|^{1/(p_{0}-1)} Sf(x)$ and
$d\mu=dx/|x|^{p_{0}/(p_{0}-1)}$, where $dx$ denotes the Lebesgue measure.
By the assumption (ii), we obtain that $T$ is bounded from
$L^{p_0}(\mathbb{R})$ to $L^{p_0}(\mathbb{R},d\mu)$. Now we shall
show that $T$ maps $L^{1}(\mathbb{R})$ to
$L^{1,\infty}(\mathbb{R},d\mu)$. Indeed, we have
$$|Tf(x)|\leq
|x|^{1/(p_{0}-1)} \|Sf\|_{\infty}\leq C|x|^{1/(p_{0}-1)} \|f\|_{1}.$$
For $\lambda>0$, a simple calculation yields
\begin{equation}\label{interpolation}
\mu\big(\{x: |Tf(x)|>\lambda\}\big)\leq \frac{C}{\lambda}\|f\|_{1}.
\end{equation}
The desired conclusion follows immediately from the Marcinkiewicz
interpolation theorem.
\end{proof}

The following lemma makes interpolation between $H_{E}^{1}\rightarrow L^{1}$ and $L^{2}\rightarrow L^{2}$
possible by the sharp function $f^{\sharp}_{E}$; see Fefferman and Stein \cite{feffermanstein} and Phong and Stein \cite{PS1986}.

\begin{lemma}\label{interpolation of sharp function}
If $F \in L^{2}$ and  $F_{E}^{\sharp}\in L^{p}$ for some $2\leq p <\infty$, then $F\in  L^{p}$ and
$$\|F\|_{p}\leq C_{p}\l\|F_{E}^{\sharp}\r\|_{p}.$$
 \end{lemma}
 \begin{proof}
 Let $f^{\sharp}$ be the well known sharp function defined by
 $$f^{\sharp}(x)=\sup_{ Q\ni x}\frac1{|Q|}\int_{Q}|f(y)-f_{Q}|dy,$$
 where $Q$ is a cube with sides parallel to the axes and $f_{Q}$ the average of $f$ over $Q$.
Observe that $f^{\sharp}(x)\le 2f^{\sharp}_{E}(x).$
By the assumption and the Fefferman-Stein
theorem (see \cite{feffermanstein} or \cite{stein}) about sharp functions, we obtain
$\|f\|_{p}\le C \|f^{\sharp}\|_{p}\leq C\|f^{\sharp}_{E}\|_{p}.$
This completes the proof.
\end{proof}

 \begin{lemma}\label{van der corput}
{\rm{(van der Corput)}} Let $I=(a,b)$ be a bounded interval on the real
line and $k\geq 1$ an integer. Suppose $\phi \in C^{k}(I)$ is real-valued
and satisfies one of the following conditions:\\
{\rm{(i)}} $k=1$, $|\phi'(t)|\geq 1$ for all $t\in I$ and $\phi'$ is monotone on $I$;\\
{\rm{(ii)}} $k\geq 2$, $|\phi^{(k)}(t)|\geq 1$ for all $t\in I$.\\
Then there exists a constant $C$, depending only on $k$ but not on $I$, such that
$$\left|\int_{I}e^{i\lambda \phi(t)}\varphi(t)dt\right|
\leq C |\lambda|^{-1/k}\left(|\varphi(b)|+\int_{I}|\varphi'(t)|dt\right)$$
for $\lambda \in \mathbb{R}$ and $\varphi \in C^{1}[a,b]$.
\end{lemma}
For the proof of this lemma and related topics, one can see \cite{stein} and \cite{carbery}.

\section{Certain Fractional Integrals}

\hspace{.5cm}We shall see that $L^p$ boundedness of $T$ in (\ref{OIO with homogeneous polynomial without cut-off}) has some connections to certain fractional integrals of Hilbert type. In particular, when $S$ is a monomial, the endpoint estimates will rely on properties of a simple class of fractional integral operators. In this section, we shall establish some mapping properties of these operators.

\begin{theorem}\label{variant of fractional integration}
Let $W_{a,b}$ be the integral operator given by
$$W_{a,b} f(x)=\int_{-\infty}^{\infty}
\big||x|_{a}^{a}-|y|_{a}^{a}\big|^{-\frac{1}b}f(y)dy$$
with $b\geq a>1$, then $W_{a,b}$ is bounded from $L^{p}$ to $L^{q}$ for $\frac{1}{p}=\frac{1}{q}+\frac{b-a}{b}$ and $1<p<\frac{b}{b-a}$.
\end{theorem}

\begin{proof}
Observe that the integral kernel is homogeneous of order $-a/b$. If $b=a>1$, then we can use the Minkowski inequality to prove the statement by a change of variables. For $b\geq a$, we obtain
  $$\|W_{a,b}f\|_{q}
  \leq \l\||f|^{1-\theta}\r\|_{b/(b-a)}\Big\|W_{a,a}\l(|f|^{b\theta/a}\r)\Big\|^{a/b}_{aq/b} \leq C
\|f\|^{1-\theta}_{b(1-\theta)/(b-a)}\|f\|_{q\theta}^{\theta},$$
with $b/a<q<\infty$. If we set $p=b(1-\theta)/(b-a)=q\theta $, i.e. $\theta=p/q$, the
desired inequality follows.
\end{proof}

Now we apply above theorem to establish $(L^{p},L^{q})$ estimates for oscillatory integral operators with polynomial phases.

\begin{theorem}\label{OIO L^p-L^q estimate (not sharp)}
\underline{} Suppose that $S$ is a real-valued homogeneous polynomial given by {\rm(\ref{homogeneous polynomial})} and $k_{min}$ as in {\rm(\ref{k_min})}. If $k_{min}\leq n/2$, then $T$  defined by {\rm (\ref{OIO with homogeneous polynomial without cut-off})} extends as a bounded operator from
$L^{2(n-k_{min})/(2n-3k_{min})}$ to $L^{2}$.
\end{theorem}
\begin{proof}
If $n=2$, $T$ reduces to the Fourier transform and hence $T$ is
bounded on $L^{2}$ by Plancherel's theorem. Now assume $n>2$. For
any $r>0$, an application of the van der Corput lemma yields
\begin{eqnarray*}
\int_{-r}^{r}|Tf(x)|^{2}dx &=& \int_{-r}^{r}\int_{\bR}\int_{\bR}
\exp\left(ia_{k_{min}}x^{n-k_{min}}\l(y^{k_{min}}-z^{k_{min}}\r)+Q(x,y,z)\right)
f(y)\overline{f(z)}dydzdx\\
&\le & C \int_{\bR}\int_{\bR}
\Big| |y|^{k_{min}}-|z|^{k_{min}}\Big|^{-1/(n-k_{min})}|f(y)f(z)|dydz\\
&\le & C \|f\|_{2(n-k_{min})/(2n-3k_{min})}^{2},
\end{eqnarray*}
where the last inequality follows from the Hardy-Littlewood-Sobolev
inequality if $k_{min}=1$ and from Theorem \ref{variant
of fractional integration} if $k_{min}>1$. Letting $r \rightarrow \infty$, we
conclude the proof.
\end{proof}

Let $S$ be a real-valued homogeneous polynomial in
$\mathbb{R}\times \mathbb{R}$ as
in (\ref{homogeneous polynomial}) with degree $n\geq 3$. Then
its partial derivative $S_{xy}''(x,y)$ can be written as
\begin{equation}\label{section 3. the Hessian of S}
S_{xy}''(x,y)=cx^{\gamma}y^{\beta}\prod_{j=1}^{m}
(y-\alpha_{j}x)^{m_{j}}\prod_{j=1} ^{s}Q_{j}(x,y),
\end{equation}
where $\alpha_{j}\neq 0$ are distinct real numbers and $Q_{j}$ are
positive definite quadratic forms. It is clear that
$\gamma=n-k_{max}-1$ and $\beta=k_{min}-1$.

In the proof of our main results, the main step is to construct an
analytic family of operators in a strip and then show these operators
satisfying suitable estimates at the boundary of the strip
considered. These operators are closely related to the Hessian
$S_{xy}''$ of $S$. As Theorem \ref{Lp boundedness with damping factor} shown, the damping factor
$|S_{xy}''|$ with suitable power gives us sharp estimates on $L^{2}$. However, it is not true generally for the endpoint estimate $H_{E}^{1}-L^{1}$. We shall see that this difference depends on whether $\beta=0$ or not in $\S$\ref{section: damped OIO}. Assume $\beta=0$. Then the following theorem gives mapping properties of the integral operator with kernels $|S''_{x,y}|^z$ with $Re(z)=-1/(n-2)$ if either (i) $\gamma>0$ and $m+s\geq 1$ or (ii) $\gamma=0$ and $m+s\geq 2$.

\begin{theorem}\label{Fraction integrals on Hardy spaces}
Suppose that $K$ is a measurable function on $\mathbb{R}^{n}\times
\mathbb{R}^{n}$ satisfying
\begin{equation}\label{K-condition-1}
\left| K\left( x,y\right) \right| \leq c\left| x\right|
^{-\theta_{0}n}\prod\limits _{k=1}^{m}\left| x-\alpha_{k}y\right|
^{-\theta_kn}
\end{equation}
 and
 \begin{equation}\label{section 3. size derivatives of K}
 \left|\nabla_{y} K\left( x,y\right)
\right| \leq C \left| x\right|^{-\theta_{0}n}\sum\limits
_{k=1}^{m}|x-\alpha_{k}y|^{-\theta_{k}n-1}\prod\limits_{j\neq
k}\left|x-\alpha_{j}y\right|^{-\theta _{j}n},
\end{equation}
where $0\leq \theta_{0}<1$, $0<\theta_{k}<1$ satisfy
$\theta_{0} +\sum _{k=1}^{m}\theta _{k}=1$ and
$\alpha_{1},\cdots,\alpha_{m}$ are distinct nonzero numbers. Let
$T_{K}$ be the integral operator given by
\begin{equation}\label{opertor-3}
T_{K}f(x)=\int_{\mathbb{R}^{n}}K(x,y)f(y)dy.
\end{equation}
Then $T_{K}$ has a bounded extension from $L^{p}$ to itself
for $1<p<\theta_{0}^{-1}$ and is also bounded from $H^{1}$ to $L^{1}$.
\end{theorem}

\begin{proof} The assumptions imply that either
$m\geq 1$ in the case $\theta_{0}>0$ or $m\geq 2$ when $\theta_{0}=0$.
It suffices to prove the theorem in the case $0<\theta_0<1$ since the treatment of
other cases is similar. To prove the $L^{p}$ boundedness of $T_K$ for
$1<p<\theta_{0}^{-1}$, we observe that
$\left| K\left( x,y\right) \right| \leq \prod _{l=1}^{n}K_{l}\left(
x_{l},y_{l}\right),$
where $x=(x_1,\cdots,x_n)$, $y=(y_1,\cdots,y_n)$ and
$ K_{l}\left( x_{l},y_{l}\right) =\left| x_{l}\right|^{-\theta
_{0}}\prod _{k=1}^{n}\left|
x_{l}-\alpha_{k}y_{l}\right|^{-\theta_{k}}.$
It follows that
\begin{equation}
\left| T_{K}f\left( x\right) \right| \leq \int_{\mathbb{R}^{n}}\prod
_{l=1}^{n}K_{l}\left( x_{l},y_{l}\right)\left| f(y) \right|dy.
\end{equation}
By Minkowski's inequality, we can reduce higher dimension $n \geq 2$
to dimension one $n=1$. For $n=1$, by a change of variables, we obtain
\begin{eqnarray*}
\left( \int_{\mathbb{R}}\left|T_{K}f\left( x\right) \right|^{p}dx\right)
^{1/p}
&\leq& \left\{\int_{\mathbb{R}}\left(\int_{\mathbb{R}}\left|
x\right|^{-\theta _{0}}\prod\limits _{k=1}^{m}\left|
x-\alpha_{k}y\right|^{-\theta_{k}}|f(y)|dy\right)^{p}dx\right\}^{1/p}\\
&\leq&  \left(\int_{\mathbb{R}}|y|^{-1/p}\prod_{k=1}^{m}\left|
1-\alpha_{k}y\right|^{-\theta_{k}}dy\right)\left\|
f\right\|_{p}.
\end{eqnarray*}
Note $0<\theta _{0}<1$, $\theta_{0} +\sum _{k=1}^{m}\theta _{k}=1$
and $1<p<\theta_{0}^{-1}$. It is easy to see that the above integral with respect to $y$ is finite.
Consequently, we have $\|T_{K}f\left( x\right) \|_p\le C\|f\|_p$ for all $1<p<\theta_{0}^{-1}$.

Now we turn to prove that $T_{K}$ is bounded from $H^{1}$ to
$L^{1}$. Suppose $a$ is any given $L^{\infty}$ atom in $H^{1}$. Then
there exists a cube $Q$ such that (i) $\supp(a)\subset Q$; (ii) $\|a\|_{\infty}\leq |Q|^{-1}$ and (iii) ${\displaystyle\int a(x)dx=0}$. It suffices to show that there exists a
 constant $C$, independent of $a$, such that $\|T_{K}(a)\|_{L^{1}}\leq
 C$. Let $d_{Q}$ be the diameter of $Q$ and $c_{Q}$ its
 center. By dilation, we assume $|\alpha_{k}|\leq 1$ for $1\leq k \leq m$. We shall now divide the proof into two cases. The first case is $|c_Q|\leq 2d_Q$. Then $|y|\leq 3d_Q$ for $y\in Q$. By the $L^p$ boundedness of $T_K$, we apply H\"{o}lder's inequality to obtain
$\left\|T_{K}\left( a\right) \right\|_{L^{1}(|x|\leq 5d_Q)}
\leq C\left|Q\right|^{1-1/p}\left\| T\left( a\right)
\right\|_{L^{p}}
\leq C$
for $1<p<\theta^{-1}$.
When $|x|>5d_Q$ and $y \in Q$, by the size assumption on $\nabla_y K$ and $|\alpha_k|\leq 1$, we obtain
\begin{eqnarray*}
|K(x,y)-K(x,c_{Q})|
&\leq&
C d_{Q}|x|^{-\theta_{0}n}
\sup_{z\in Q}\left(\sum_{k=1}^{n}|x-\alpha_k z|^{-\theta_k n-1}
\prod_{j\neq k}|x-\alpha_{j}z|^{-\theta_{j}n}\right)\\
&\leq& Cd_Q|x|^{-n-1}.
\end{eqnarray*}
Hence it follows that
\begin{equation*}
\int_{|x|>5d_Q}|K(x,y)-K(x,c_Q)|dx
\leq Cd_Q \int_{|x|> 5d_Q}|x|^{-n-1}dx< \infty.
\end{equation*}
Thus we have
\begin{eqnarray*}
\|T_{K}a\|_{L^1(|x|>5d_Q)}&=&\int_{|x|>5d_Q}\left|\int_{Q}\Big(K(x,y)-K(x,c_{Q})\Big)a(y)dy\right|dx\\
 &\leq &\left(\sup_{y \in Q}\int_{|x|>5d_Q}|K(x,y)-K(x,c_{Q})|dx\right)\int_{Q}|a(y)|dy \\
 &\leq& C.
\end{eqnarray*}

The second case is $|c_Q|> 2d_Q$. It is easy to see that
$|c_Q|/2\leq |y|\leq 3|c_Q|/2$ for all $y\in Q$. This observation implies
\begin{equation*}
\sup_{y\in Q}\int_{|x|\leq 5|c_Q|}\left|K(x,y)\right|dx\leq C<\infty
\end{equation*}
by the assumption $0<|\alpha_k|\leq 1$. Thus the integral of
$|T_K(a)|$ over the ball $|x|\leq 5|c_Q|$ is bounded by a constant $C$.
For $|x|\geq 5|c_Q|$, $|K(x,y)-K(x,c_Q)|$ has the same upper bound
$Cd_Q|x|^{-n-1}$ uniformly for $y\in Q$. Hence the integral
$|T_K(a)|$ over $|x|\geq 5|c_Q|$ is also less than a bound independent
of $a$.

Combing above estimates, we conclude the proof of the theorem.
\end{proof}

\section{Oscillatory Integral Operators On $H_E^1$}

\hspace{.5cm} Let $P$ be a real-valued polynomial on
$\mathbb{R}^{n}\times\mathbb{R}^{n}$. Let $T_{P}$ be an oscillatory
integral operator as in (\ref{osc singular integral with polynomial phases}) with the kernel
$K$ as in Theorem \ref{Fraction integrals on Hardy spaces}. In this
section, our purpose is to establish the boundedness of $T_P$ from
$H_E^1$ to $L^1$. This result will serve as an endpoint estimate for
an analytic family of operators in the proof of Theorem
\ref{main theorem 1}. If $P$ can be written as $P(x,y)=P_1(x)+P_2(y)$,
it is easy to see that $T_P$ is bounded from $H_E^1$ to $L^1$. In fact, $H_{E}^{1}$ is isomorphic to $H^1$ by a multiplication of $\exp{(-iP_{2}(y))}$.

\begin{theorem}\label{osc fractional operators with polynomial phases}
Suppose that $P$ is a real-valued polynomial and that
$T_{P}$ is given by {\rm(\ref{osc singular integral with
polynomial phases})} with $K$ as in Theorem {\rm \ref{Fraction
integrals on Hardy spaces}}. Then $T_{P}$ is a bounded operator from
$L^{p}$ to itself for $1<p<\theta_{0}^{-1}$, and  $T_{P}$ has a
bounded extension from $H^{1}_{E}(\mathbb{R}^{n};P)$ to
$L^{1}(\mathbb{R}^{n})$. Moreover, the operator norms of $T_{P}$
depend only on the degree of $P$ but not on its coefficients.
\end{theorem}

We first introduce a useful lemma which was essentially proved by Ricci
and Stein in \cite{ricci}. Some related topics were systematically
studied in \cite{carbery}.

\begin{lemma}\label{deacy of OIO with P phase}
Let $T_{\lambda}$ be the oscillatory integral operator defined by
\begin{equation*}
T_{\lambda}f(x)=\int_{\mathbb{R}^{n}}e^{i\lambda
P(x,y)}\phi(x,y)f(y)dy,
\end{equation*}
where $P$ is a real-valued polynomial with degree $d\geq 2$, $\phi
\in C_{0}^{\infty}$ and $\lambda \in \mathbb{R}$. If there exist
multi-indices $\beta_1$ and $\beta_2$ satisfying
\begin{equation}\label{section 4. mixed par-deri geq 1}
\left|\frac{\partial^{d}P}{\partial^{\beta_{1}}x
\partial^{\beta_{2}}y}\right|\geq 1,
\end{equation}
with $0<|\beta_1|<d$ and $|\beta_1|+|\beta_2|=d$, then there exists
some $\delta>0$, depending only on the degree $d$ of $P$, such that
the following decay estimate holds,
\begin{equation}
\|T_{\lambda}f\|_{2}\leq C|\lambda|^{-\delta}\|f\|_{2}.
\end{equation}
\end{lemma}

\begin{proof}
For $1<p<\theta_{0}^{-1}$, the $L^{p}$ boundedness of $T_{P}$  follows
immediately from Theorem \ref{Fraction integrals
on Hardy spaces}, since we can take absolute value in the integral and the resulting operator is bounded on $L^p(\bR^n)$. Now we turn our attention to prove the boundedness of $T_P$ from $H_E^1$
to $L^1$ by induction on the degree $d$ of $P$. For $d=0$, this statement
is just Theorem \ref{Fraction integrals on Hardy spaces}. For $d=1$, $P(\cdot,\cdot)$ is degenerate in the sense that it can be decomposed as the sum of two polynomials $P_{1}(x)$ and $P_{2}(y)$. The statement is also true.

Assume the theorem is true for all polynomials of degree $l$ not greater
than $d-1$. We shall prove that it is also true for $l=d$. Write
\begin{equation*}
P(x,y)
=\sum\limits_{|\beta_1|+|\beta_2|=d}c_{\beta_1,\beta_2}
x^{\beta_{1}}y^{\beta_2}+Q(x,y)=P_{d}(x,y)+Q(x,y)
\end{equation*}
with $deg(Q)\leq d-1$. Assume that $P_{d}$ is not degenerate without loss
of generality. For the same reason, we may assume that the coefficients
of the pure $x^{\alpha}$ and $y^{\beta}$ terms with $|\alpha|=|\beta|=d$
are zero. Given any $H^{1}_{E}$ atom $a$ associated to $P$, the aim now
is to show that $\|T_{P}a\|_{1}\leq C$ for some $C<\infty$ independent of
$a$. Let $Q$ be the cube associated with $a$, with the center $c_{Q}$
and diameter $d_{Q}$. For convenience, we divide the proof into two
cases.

Case I. $|c_Q|\leq 2d_Q$.

Set $N=10+\max_{k}{|\alpha_{k}|}$. Since $T_P$ is bounded on
$L^{p}$ for $1 < p < \theta _{0}^{-1}$, we obtain
\begin{equation*}
\left\| T a\right\|_{L^{1}\left(|x|\leq 5Nd_Q\right)}
\leq
C\left|Q\right|^{1-1/p}\left\| T a\right\|_{L^{p}}
\leq C.
\end{equation*}
For $|x|>5Nd_Q$, we may write
\begin{eqnarray*}
T_{P}(a)(x)&=&\int_{Q}e^{iP(x,y)}K(x,y)a(y)dy\\
         &=&\int_{Q}e^{iP(x,y)}\Big(K(x,y)-K(x,c_Q)\Big)a(y)dy
         +K(x,c_Q)\int_{Q}e^{iP(x,y)}a(y)dy\\
         &=& I_1 +I_2.
\end{eqnarray*}
Recall we have proved that in Theorem
\ref{Fraction integrals on Hardy spaces}
\begin{equation}\label{variant of Hormander condition}
\sup\limits_{y \in
Q}\int_{|x|>5Nd_Q}\left|K(x,y)-K(x,c_{Q})\right|dx \leq C<\infty
\end{equation}
which implies
\begin{equation*}
\int_{|x|>5Nd_Q}\left|I_{1}\right|dx \leq C\int_{Q}|a(y)|dy\leq
C<\infty.
\end{equation*}
Before applying Lemma \ref{deacy of OIO with P phase} to estimate
$I_2$, we need an additional argument to show that the bounds are
independent of the coefficients of $P$. Recall that the coefficients
of the pure $x^\alpha$ and $y^\beta$ terms are assumed to be zero
for all multi-indices $|\alpha|=|\beta|=d$. Decompose $P$ as
\begin{equation*}\label{decomposition of the polynomial}
P(x,y)=\sum\limits_{|\alpha|+|\beta|=d}
         c_{\alpha,\beta}(x-c_Q)^{\alpha}(y-c_Q)^{\beta}+R(x,y),
\end{equation*}
where the degree of $R$ is less than $d$. Then it is easy to see
$P(c_{Q},y)=R(c_Q,y)$. By the induction hypothesis, $\|T_{R}(a)\|_{1}
\leq C$. For $|x|>5Nd_Q$, we can also write $T_R(a)$ as
\begin{equation*}
T_{R}(a)(x)=\int_{Q}e^{iR(x,y)}\Big(K(x,y)-K(x,c_Q)\Big)a(y)dy
         +K(x,c_Q)\int_{Q}e^{iR(x,y)}a(y)dy
\end{equation*}
By (\ref{variant of Hormander condition}), it follows that
\begin{equation}\label{auxilliary estimate of OIO on Hardy space}
\int_{|x|>5Nd_Q}\left|K(x,c_{Q})\int_{Q}e^{iR(x,y)}a(y)dy\right|dx
\leq C,
\end{equation}
with the bound $C$ independent of the coefficients of $R$. For
$|x|>5Nd_Q$, then
$|x-c_Q|\geq 4Nd_Q$. It is convenient to observe that
(\ref{auxilliary estimate of OIO on Hardy space}) is still true
with $|x|>5Nd_Q$ replaced by $|x-c_Q|\geq 4N d_Q$. For $t>0$, write
\begin{eqnarray*}
& &\int_{4Nd_{Q}\leq |x-c_{Q}|\leq t}
\left|K(x,c_Q)\int_{Q}e^{iP(x,y)}a(y)dy\right|dx\\
 &&=
\int_{4Nd_{Q}\leq |x-c_{Q}|\leq t}\left|K(x,c_Q
)\int_{Q}\l(e^{iP(x,y)}-e^{iR(x,y)}\r)a(y)dy\right|
dx\\
&&\qquad + \int_{4Nd_{Q}\leq |x-c_{Q}|\leq t}
\left|K(x,c_Q)\int_{Q}e^{iR(x,y)}a(y)dy\right|dx\\
 &&=
I_{2,1}+I_{2,2}.
\end{eqnarray*}
It is clear that $|I_{2,2}|\leq C$ since (\ref{auxilliary estimate of
OIO on Hardy space}) is still true with $|x|>5Nd_Q$ replaced by $|x-c_Q|\geq 4N d_Q$.
Note that $|x-\alpha_kc_{Q}|\approx |x-c_{Q}|$ for
$|x-c_Q|>4Nd_Q$, we have
\begin{eqnarray*}
|I_{2,1}|&\leq&
\sum\limits_{\alpha,\beta}|c_{\alpha,\beta}|
\|a\|_{L^{1}}\int_{C_{1}d_{Q}\leq |x-c_{Q}|\leq
t}|x-c_{Q}|^{|\alpha|}|y-c_{Q}|^{|\beta|}|K(x,c_Q)|dx\\
&\leq&
\sum\limits_{\alpha,\beta}|c_{\alpha,\beta}|
d_{Q}^{d-|\alpha|}\int_{C_{1}d_{Q}\leq |x-c_{Q}|\leq
t}|x-c_{Q}|^{|\alpha|-n}dx\\
&\leq&
\sum\limits_{\alpha,\beta}|c_{\alpha,\beta}|
t^{|\alpha|}d_{Q}^{d-|\alpha|},
\end{eqnarray*}
where the summations are taken over all multi-indices $\alpha$ and $\beta$ satisfying $|\alpha|+|\beta|=d$ with $0<|\alpha|<d$.
Take
$$t^{-1}=\l(\max|c_{\alpha,\beta}|
d_{Q}^{d-|\alpha|}\r)^{1/|\alpha|}>0,$$
where the maximum is taken over all multi-indices appearing in the above
summations. It follows that $|I_{2,1}|\leq C$ with $C$
depending only on the degree $d$ and the dimension $n$. It remains
to show that
\begin{equation}\label{major estimate of OIO on Hardy space}
\int_{|x-c_{Q}|\geq \max\{t,4Nd_{Q}\}}
\left|K(x,c_Q)\int_{Q}e^{iP(x,y)}a(y)dy\right|dx\leq C.
\end{equation}
We have pointed out that $|x-c_{Q}|\geq 4N d_{Q}$ implies
$|K(x,c_{Q})|\approx |x-c_{Q}|^{-n}.$
By this observation and Lemma
\ref{deacy of OIO with P phase}, we can assume now that $Q$ is
the unit cube centered at the origin. Indeed, we have
\begin{eqnarray*}
& &\int_{|x-c_{Q}|\geq \max\{t,4N d_{Q}\}}
\left|K(x,c_Q)\int_{Q}e^{iP(x,y)}a(y)dy\right|dx\\
&&\leq \sum_{j=0}^{\infty}\int_{{2^{j-1}t\leq |x-c_Q|<2^{j}t}}|x-c_{Q}|^{-n}
\left|\int_{Q}e^{iP(x,y)}a(y)dy\right|dx.
\end{eqnarray*}
Let $Q_0$ be the unit cube centered at the origin and $l$ the side
length of $Q$. By a change of variables, we obtain
\begin{eqnarray*}
&&\int_{2^{j-1}t\leq |x|< 2^{j}t}
|x|^{-n}
\left|\int_{Q_{0}}e^{iP(x+c_{Q},ly+c_{Q})}l^{n}a(ly+c_{Q})dy\right|dx\\
&\leq & C \int_{|x|\leq
2}\left|\int_{Q_{0}}e^{iP(2^{j}tx+c_{Q},ly+c_{Q})}l^{n}a(ly+c_{Q})dy
\right|dx\\
&=& C\big\|U_{j}\l(l^{n}a(l(\cdot)+c_{Q})\r)\big\|_{L^{2}(|x|\leq
2)}
\end{eqnarray*}
where $U_j$ is the operator given by
$$U_jf(x)=\int_{Q_{0}}e^{iP(2^{j}tx+c_{Q},ly+c_{Q})}f(y)dy.$$
By the choice of $t$, there exist multi-indices $\alpha$ and $\beta$
satisfying $0<|\alpha|,|\beta|<d$ and $|\alpha|+|\beta|=d$ such that
$$\left|\frac{\partial^d }{\partial x^\alpha \partial y^\beta}
P(2^{j}tx+c_{Q},ly+c_{Q})\right|\geq C 2^{j|\alpha|}$$
with the constant $C$ depending only on $n$ and $\alpha$.
By Lemma \ref{deacy of OIO with P phase}, the power decay property of
$\|U_j\|_{L^2\rightarrow L^2}\leq C 2^{-j\delta}$ for some $\delta>0$ implies the desired
estimate (\ref{major estimate of OIO on Hardy space}).\\

 Case II. $|c_{Q}|> 2 d_{Q}$.\\
By the size condition (\ref{K-condition-1}) imposed on the kernel $K$,
it follows that
\begin{equation*}
\sup_{y \in Q}\int_{|x|\leq 2N |c_{Q}|}\left|K(x,y)\right|dx\leq C
<\infty
\end{equation*}
with the constant $C$ independent of $Q$. Let
$t^{-1}=\l(\max|c_{\alpha,\beta}|d_{Q}^{d-|\alpha|}\r)^{1/|\alpha|}>0$
as above. By a similar argument, it is true that
\begin{equation*}
\int_{2N |c_{Q}|<|x|<t}\left|\int_{Q}e^{iP(x,y)}K(x,y)a(y)dy\right|dx
\leq C<\infty.
\end{equation*}
At the same time, the integral of $|T_P(a)|$ over $|x|>2N |c_Q|$ is not
greater than
\begin{eqnarray*}
&&\int_{|x|>2N |c_{Q}|}\left|\int_{Q}e^{iP(x,y)}[K(x,y)-K(x,c_{Q})]a(y)dy\right|dx\\
&&+
\int_{|x|>2N |c_{Q}|}|K(x,c_{Q})|\left|\int_{Q}e^{iP(x,y)}a(y)dy\right|dx\\
&=& J_{1}+J_{2}.
\end{eqnarray*}
To estimate $J_{1}$, first observe that the assumption
(\ref{section 3. size derivatives of K}) of the kernel $K$ implies
$$|K(x,y)-K(x,c_{Q})|\leq Cd_{Q}|x|^{-n-1}$$
for all $|x|>2N|c_Q|$. Hence we have $J_{1}\leq C$. Now it remains to show that
\begin{equation}
\int_{|x|>\max\{t,2N |c_{Q}|\}}|K(x,c_{Q})|\left|\int_{Q}e^{iP(x,y)}a(y)dy\right|dx
\leq C<\infty.
\end{equation}
The proof is similar as above. Hence the theorem is completely proved.
\end{proof}

The following theorem deals with a class of translation invariant operators.

\begin{theorem}\label{H1 boundedness of trans-inv OIO}
 Let
\begin{equation}\label{4-30-operator}
Uf(x)=\int_{-\infty}^{\infty}e^{i(x-y)^n}(1+|x-y|)^{z}f(y)dy,
\end{equation}
where $n\geq 2$ is an integer and  $\Re(z)=-1$.
Then $U$ is bounded on $L^{2}$ and maps $H^1$ boundedly into $L^1$ with
the bounds less than a constant multiple of $(1+|\Im(z)|)^{2}$.
\end{theorem}
As mentioned in the introduction, the theorem is essentially contained in
\cite{jurkatsampson} (at least for even $n$). The argument in
\cite{jurkatsampson} is also applicable here. We shall give a new proof
of the $L^2$ boundedness at the end of $\S$5. By the argument in the proof of above theorem, we can show that $U$ is bounded from $H^1$ to $L^1$.

\section{Damped Oscillatory Integral Operators}\label{section: damped OIO}

\hspace{.5cm}In $\S$4, we have considered a class of oscillatory integral operators with critical negative power and obtained the endpoint estimate from $H_{E}^{1}$ to $L^{1}$. To put $T$
defined by (\ref{OIO with homogeneous polynomial without cut-off}) into a family of analytic operators, we shall insert a damped factor into $T$. For this purpose, we shall study operators of the following form
\begin{equation}\label{OIO with damped factors}
Wf(x)=\int_{-\infty}^{\infty}e^{iS(x,y)}|D(x,y)|^{z}f(y)dy
\end{equation}
with suitably chosen $z$ and the damped factor $D$ being determined by $S''_{xy}$.

In this section, our main result is Theorem \ref{OIO with damped factor on L2}. It serves as an endpoint estimate of the operator $W$. The other endpoint estimate for $W$ has been obtained in $\S$4 except some special cases. More precisely, when the Hessian is of form $S''_{xy}=cy^{\beta}(y-\alpha x)^{n-2-\beta}$,
the treatment is different and we shall deal with this case separately. For related results about damped oscillatory integral operators, we refer the reader to \cite{PS1998} and \cite{pyang}. A general class of weighted oscillatory integral operators had been studied in \cite{pyang}. The region obtained in \cite{pyang} is given by an infinite intersection whose boundary is obscure. The damped factor $D$ and its critical damping exponent are explicitly given in this section.

\begin{theorem}\label{OIO with damped factor on L2}
Assume the Hessian $S_{xy}^{''}$ as in {\rm(\ref{section 3. the Hessian of S})}. Let $W$ be given by {\rm(\ref{OIO with damped factors})} with
\begin{equation}\label{section 5. standard damping factor}
D(x,y)=x^\gamma\prod^m_{j=1}(y-\alpha_jx)^{m_j}
\prod^s_{j=1}Q_j(x,y)
\end{equation}
and
$$\Re(z)=a_{\beta}=\frac1{2(\beta+1)}\frac{n-2(\beta+1)}{n-\beta-2}.$$
Then $W$ extends as a bounded operator on $L^2(\bR)$ with the operator norm $\|W\|$ less than a constant multiple of $(1+|\Im(z)|)^{2}$.
\end{theorem}

\noindent{\bf Remark.}  For $\beta=0$, this result is the same as Theorem \ref{Lp boundedness with damping factor}. When $\beta=(n-2)/2$ and $z=0$, the statement is contained in Theorem
\ref{OIO L^p-L^q estimate (not sharp)} since $k_{min}=n/2$ and $a_{\beta}=0$. If $\beta=n-2$, $a_{\beta}$ is not well defined. In this case, we put the damping factor $D\equiv1$ and $|D|^z=1$. By insertion of a smooth cut-off, we obtain
\begin{equation}\label{OIO with damped factor and cut-off}
W_\lambda
f(x)=\int_{\bR}e^{i\lambda S(x,y)}|D(x,y)|^{z}\varphi(x,y)f(y)dy
\end{equation}
for $\lambda\in \bR$ and $\varphi\in C^\infty_0(\bR\times\bR)$. By dilation, the $L^2$ boundedness is equivalent to the following decay estimate
\begin{equation}\label{decay estimate of OIO with Damped factor}
\|W_\lambda f\|_{L^2}\le
C|\lambda|^{-\frac1{2(1+\beta)}}\|f\|_{L^2}.
\end{equation}
Our proof uses the techniques introduced in \cite{PS1994}, \cite{PS1997} and \cite{PS1998}. By
a decomposition of the singular variety $\{S''_{xy}(x,y)=0\}$, we can write $W_\lambda=\sum_{k,l\ge0}W_{k,l}$. By the almost orthogonality principle, we shall balance size and oscillatory estimates to obtained the sharp $L^2$ decay rate. The oscillatory estimates rely on the following operator version of van der Corput lemma.

\begin{lemma}\label{Operator version of Van der Corput lemma}
{\rm(\cite{PS1997})} Assume $T_\lambda$ and $S$ as in the introduction.
Suppose that $\phi(x)$ is monotone and continuous
 on the interval $[\alpha,\beta]$. For some $\delta>0$, the cut-off
 function $\varphi$ and the phase $S$
 satisfy\\

{\rm(i)} $\supp(\varphi)\subset \big\{(x,y):\;\phi(x)\leq y \leq \phi(x)+\delta, \;\alpha \leq x \leq \beta\big\}$;\\

{\rm(ii)} $\displaystyle{\l|\partial_{y}^k \varphi(x,y)\r|\leq C\delta^{-k}}$ for $k=0,1,2$;\\

{\rm(iii)} $\displaystyle{\mu \leq \big|S_{xy}^{''}\big|\leq C\mu}$ on
the curved box
$\{(x,y):\;\phi(x)\leq y \leq \phi(x)+\delta, \;\alpha \leq x \leq \beta\}$\\
for some $\mu >0$. Then
$$\|T_\lambda f\|_{L^2}\le
C|\lambda \mu|^{-1/2}\|f\|_{L^2}$$
with the bound independent of $f$.
\end{lemma}

We also need a lemma which measures the orthogonality in the decomposition of the operator $W_{\lambda}$. Results of this type is of fundamental importance in the study of damped
oscillatory integral operators; see \cite{PS1998}. The following lemma is a simplified version of those given in \cite{PS1998}.

Let $T_1$ and $T_2$ be two oscillatory integral operators given by
\begin{equation}
T_{j}f(x)=\int_{-\infty}^{\infty}e^{i\lambda S(x,y)}\varphi_{j}(x,y)
f(y)dy,\quad j=1,2
\end{equation}
where $S$ is a real-valued homogeneous polynomial given as in
(\ref{homogeneous polynomial}) and $\varphi_{j}\in C_0^\infty$. Suppose
$\varphi_j$ is supported in a parallelogram which has two sides parallel
to the $y-axis$ with length $\delta_j$ for each $j$. Denote these
two parallelograms by $\Omega_1$ and $\Omega_2$. Then
$$\Omega_{j}
=\{(x,y):a_j\leq x\leq b_j,\quad l_j(x)\leq y \leq l_j(x)+\delta_j\}$$
where $y=l_j(x)$ is a line in the plane.
\begin{lemma}\label{section 5. orthogonality structure}
Let $T_1$ and $T_2$ be defined as above. Let $\Omega_1^{\ast}$ be
the expanded parallelogram
$$\Omega_{1}^{\ast}(c)
=\{(x,y):a_1\leq x\leq b_1,
\quad l_1(x)-c\delta_{1}\leq y \leq l_1(x)+(1+c)\delta_1\}$$
with $c>0$. If there exist positive numbers $c_1$ and $c_2$ with
$c_1<c_2$ such that $\Omega_2\subset \Omega_1^{\ast}(c_2)$, the Hessian
$S''_{xy}$ does not change sign in $\Omega_1^{\ast}(c_2)$ and
$$\min_{\Omega_1^{\ast}(c_1)}|S''_{xy}|\geq \mu >0
\quad and \quad \max_{\Omega_1^{\ast}(c_2)}|S''_{xy}|\leq C\mu,$$
then we have
\begin{equation*}
\|T_1T_2^{\ast}\|_{L^2\rightarrow L^2}\leq C(|\lambda|\mu)^{-1}\prod_{j=1}^{2}\left(
\max_{\Omega_{j}}\sum_{k=0}^{2}\delta_j^k|\partial_{y}^k\varphi_j(x,y)|\right).
\end{equation*}
\end{lemma}
\noindent{\bf Remark.} If $\Omega_1$ and $\Omega_2$ are two rectangles
with sides parallel to the axes, the above assumptions can be relaxed.
The assumptions in the lemma can be replaced by the following
conditions:\\
(i) Let $\Omega_1^\ast(c)$ be the expanded rectangle only in the
$x$-dimension but with the $y$-dimension unchanged,
$$\Omega_1^\ast(c)=\{(x,y):a_1-c(b_1-a_1)\leq x
\leq b_1+c(b_1-a_1),\quad \alpha_1\leq y \leq \beta_1\}$$
with $c>0$. Then assume that
$$\min_{\Omega_{1}^\ast(c_1)}|S''_{xy}|\geq \mu >0,
\qquad \Omega_2\subset \Omega_1^{\ast}(c_2)$$
for some $c_1>0$ and $c_2>0$;\\
(ii) Let $\mathcal{R}$ be the rectangle consisting of all segments
joining two points $(x_1,y)\in \Omega_1$ and $(x_2,y)\in \Omega_2$.
The Hessian $S''_{xy}$ does not change sign on $\mathcal{R}$ and
$\max_{\mathcal{R}}|S''_{xy}|\leq C\mu$.

The proof of above two lemmas relies on a basic property of polynomials.
More precisely, for an arbitrary given polynomial $P$ in $\mathbb{R}$,
we have
$$\sup_{x\in I^\ast}|P^{(k)}(x)|
\leq C|I|^{-k}\sup_{x\in I}|P(x)|$$
where the bound $C$ depends only on the degree of $P$ but not on
the choice of the interval $I$ and $I^\ast$ is the interval concentric
with $I$ but dilated by the factor $2$. A general concept of
polynomial-like functions was introduced in \cite{PS1997}.

We also use Schur's lemma frequently. For convenience, we state it as
follows.
\begin{lemma}\label{Schur's Lemma}
Suppose that $K(\cdot,\cdot)$ is measurable in $\mathbb{R} \times
\mathbb{R}$ satisfying
$$
\sup_{x}\int|K(x,y)|dy \leq A
\quad \textrm{and} \quad
 \sup_{y}\int|K(x,y)|dx \leq B
$$
for some $A,B>0$. Then the integral operator with kernel $K$ is
bounded on $L^{2}$ with bound not greater than $\sqrt{AB}$.
\end{lemma}

With above preliminaries, we turn to the proof of Theorem
\ref{OIO with damped factor on L2}. The symbol $\|V\|$ means the
norm of the operator $V$ on $L^2$. We also use $A\approx B$ to
mean that $C_1A\leq B \leq C_2A$ for some constants $C_1,C_2>0$.\\

\begin{proof}
 Choose a smooth function $\Phi$ satisfying
$\supp\Phi\subset[1/2,2]$ and
$\sum_{j\in\bZ}\Phi(x/{2^j})=1$ for all $x>0$.
If the support of $\varphi$ is sufficiently small, then we can divide
$W_\lambda$ in (\ref{OIO with damped factor and cut-off}) into four
parts of the form
$$\int_{\bR}e^{i\lambda S(x,y)}|D(x,y)|^{z}\chi_{\{\pm x>0\}}(x)\chi_{\{\pm y>0\}}(y)\varphi(x,y)f(y)dy,$$
where $\chi_A$ is the characteristic function of the set $A$.

We shall prove that each of above four operators satisfies
the desired estimate (\ref{decay estimate of OIO with Damped
factor}). Since the argument is similar, it suffices to show that the
desired estimate holds for one of these operators. For the estimate of
$W_\lambda$ in the first quadrant, i.e., $x>0$, $y>0$,
decompose $W_\lambda$ as
$$W_\lambda f(x)=\sum_{k,l \geq 0}W_{k,l}f(x),$$
where
\begin{equation}\label{section 5. def of W(k,l)}
 W_{k,l}f(x)=\int_{\bR}e^{i\lambda
S(x,y)}|D(x,y)|^{z}\Phi_k(x)\Phi_l(y)\varphi(x,y)f(y)dy
\end{equation}
with $\Re(z)=a_\beta$ and $\Phi_k(x)=\Phi(2^k x)$. We shall
use $W_{k,l}^{++}$ to denote the decomposition in the first quadrant.
Likewise, $W_{k,l}^{\sigma_1,\sigma_2}$ is defined similarly with
$\sigma_1,\sigma_2=\pm$.

{\bf Range $k\geq l+N_0$.}\\
Assume $N_0>0$ is a large number such that
$x\approx2^{-k}, y\approx2^{-l}$ imply
$|y-\alpha_ix|\approx2^{-l}$ if $k\ge l+N_0$. In this case, the proof is
somewhat different depending on whether $\gamma$ equals 0 or not. We
begin with the simpler case $\gamma=0$. It is convenient to estimate
the summation of operators $W_{k,l}^{++}+W_{k,l}^{-+}$ in the range
$k\geq l+N_0$. For each $l\geq 0$, let
$U_l=\sum_{k\geq l+N_0}(W_{k,l}^{++}+W_{k,l}^{-+})$. By the support of $\Phi$, we see that
$U_lU_{l'}^{\ast}=0$ unless $|l-l'|\leq 1$. For $U_l^{\ast}U_{l'}$, we
claim that $\|U_l^{\ast}U_{l'}\|$ is bounded by
$C_z|\lambda|^{-\frac{1}{\beta+1}}2^{-|l-l'|\delta}$ with some
$\delta>0$. Observe that $U_l$ is supported in the rectangle
$|x|\leq 2^{-N_0+1}2^{-l}$ and $y\approx 2^{-l}$. Without loss of
generality, we may assume $l\leq l'$. It is easy to verify that the
assumptions in the remark after
Lemma \ref{section 5. orthogonality structure} are satisfied with the
roles $x$ and $y$ reversed. Hence we have
\begin{equation*}
\|U_l^{\ast}U_{l'}\|
\leq A= C(1+|\Im(z)|)^4(|\lambda|2^{-l(n-2)})^{-1}
2^{-l'(n-2-\beta)a_{\beta}}2^{-l(n-2-\beta)a_{\beta}}.
\end{equation*}
By Schur's lemma, it follows that
\begin{equation*}
\|U_l^{\ast}U_{l'}\|
\leq B= C2^{-l}2^{-l'}2^{-l(n-2-\beta)a_{\beta}}
2^{-l'(n-2-\beta)a_{\beta}}.
\end{equation*}
Combing above two estimates, we see that $\|U_l^{\ast}U_{l'}\|$ is not greater than $A^{\theta}B^{1-\theta}$ for all $0\leq \theta \leq 1$. Taking $\theta=1/(\beta+1)$, we obtain
\begin{equation*}
\|U_l^{\ast}U_{l'}\|\leq C(1+|\Im(z)|)^{4\theta}
|\lambda|^{-\frac{1}{\beta+1}}2^{-|l-l'|(n-2)/[2(\beta+1)]}.
\end{equation*}
By the almost orthogonality lemma, we see that $\|\sum_lU_l\|$ is not
greater than a constant multiple of
$(1+|\Im(z)|)^2|\lambda|^{-\frac{1}{2(\beta+1)}}$.

Now assume $\gamma>0$. Denote by $k\wedge k'$ the minimum of $k$ and
$k'$. Then $W_{k,l}W_{k',l'}^{\ast}$ with $|l-l'|\leq 1$ satisfies the
following estimates,
\begin{eqnarray*}
&&\|W_{k,l}W_{k',l'}^{\ast}\|
\leq C(1+|\Im(z)|)^4
\left(|\lambda|2^{-(k\wedge k')\gamma}2^{-l(n-2-\gamma)}\right)^{-1}
2^{-(k+k')\gamma a_{\beta}}2^{-(l+l')(n-2-\beta-\gamma)a_{\beta}}\\
&&\|W_{k,l}W_{k',l'}^{\ast}\|
\leq
C 2^{-(k+k')/2}2^{-(l+l')/2}
2^{-(k+k')\gamma a_{\beta}}2^{-(l+l')(n-2-\beta-\gamma)a_{\beta}}.
\end{eqnarray*}
A convex combination of above estimates gives
\begin{equation*}
\|W_{k,l}W_{k',l'}^{\ast}\|\leq C(1+|\Im(z)|)^{4/(\beta+1)}
|\lambda|^{-\frac{1}{\beta+1}}2^{-|k-k'|\delta}
\end{equation*}
with $\delta=\gamma a_{\beta}+\beta/[2(\beta+1)]>0$. Indeed, it is easy to verify that
$$\delta=\frac{1}{2(\beta+1)}\frac{(n-2)(\beta+\gamma)-(\beta+\gamma)^2+\gamma^2}{n-2-\beta}>0.$$
 Similarly, by reversing $x$ and $y$ in Lemma \ref{section 5. orthogonality structure}, it is also
true that for $|k-k'|\leq 1$
\begin{equation*}
\|W_{k,l}^{\ast}W_{k',l'}\|\leq C(1+|\Im(z)|)^{4/(\beta+1)}
|\lambda|^{-\frac{1}{\beta+1}}
2^{-|l-l'|\delta}
\end{equation*}
with $\delta=(n+\beta-2\gamma-2)/[2(\beta+1)]+\gamma a_{\beta}$. We see
$\delta>0$ unless $\gamma=n-2$. In the case $\gamma=n-2$, $\delta=0$ and we shall treat it separately. We shall rather estimate $W_{k}=\sum_{l=-\infty}^{\infty} (W_{k,l}^{++}+W_{k,l}^{+-})$.
Also $W_k^{\ast}W_{k'}=0$ if $|k-k'|\geq 2$. By the oscillatory estimate, we obtain
\begin{equation*}
\|W_{k}W_{k'}^{\ast}\|\leq C(1+|\Im(z)|)^4|\lambda|^{-1}2^{-|k-k'|(n-2)/2}.
\end{equation*}

{\bf Range $l\geq k+N_0$.}\\
Now we turn to the estimate of $\|\sum W_{k,l}\|$ with summation being taken over $k\leq l-N_0$ for a sufficiently large $N_0>0$. For arbitrary two pairs $(k,l)$ and $(k',l')$ in the range, we also have,
$|l-l'|\leq 1$,
\begin{eqnarray*}
&&\|W_{k,l}W_{k',l'}^{\ast}\|
\leq
C2^{-(k+k')/2}2^{-(l+l')/2}
2^{-(k+k')(n-2-\beta)a_{\beta}}\\
&&\|W_{k,l}W_{k',l'}^{\ast}\|
\leq
C(1+|\Im(z)|)^4
\left(|\lambda|2^{-l\beta}2^{-(k\wedge k')(n-2-\beta)}\right)^{-1}
2^{-(k+k')(n-2-\beta)a_{\beta}}.
\end{eqnarray*}
By a convex combination of the above two inequalities, we obtain
$$\|W_{k,l}W_{k',l'}^{\ast}\|
\leq
C(1+|\Im(z)|)^{4/(\beta+1)}|
\lambda|^{-\frac{1}{\beta+1}}
2^{-|k-k'|\delta}$$
with $\delta=(n-\beta-2)/[2(\beta+1)]$. By the same argument, we also have
$$\|W_{k,l}^{\ast}W_{k',l'}\|
\leq
C(1+|\Im(z)|)^{4/(\beta+1)}|
\lambda|^{-\frac{1}{\beta+1}}
2^{-|l-l'|\beta/[2(\beta+1)]},\quad |k-k'|\leq 1.$$
For the special case $\beta=0$, we
rather estimate $V_k=\sum_{l}(W_{k,l}^{++}+W_{k,l}^{+-})$ over the
range $l\geq k+N_0$. The assumptions in the remark of Lemma
\ref{section 5. orthogonality structure} are satisfied by $V_k$ and
$V_{k'}$ and then the resulting estimate is given by
$\|V_kV_{k'}^{\ast}\|\leq
C(1+|\Im(z)|)^{4/(\beta+1)}|
\lambda|^{-\frac{1}{\beta+1}}
2^{-|k-k'|(n-2)/2}.$
Observe that $V_k^{\ast}V_{k'}=0$ for $|k-k'|\geq 2$. Hence
$\|\sum V_k\|\leq C(1+|\Im(z)|)^2|\lambda|^{-\frac{1}{2(\beta+1)}}$ by the almost orthogonality principle.

{\bf Range $|k-l|\leq N_0$.}\\
 For $|k-l|\le N_0$, we shall further decompose $W_{k,l}$ as
$$W_{k,l}=\sum\limits_{d_\ge0}W_{k,l,d}
=\sum\limits_{d}\int_{\bR}e^{i\lambda S(x,y)}|D(x,y)|^{z}
\Phi_k(x)\Phi_l(y)\Phi_{d}(y-\alpha_1x)\varphi(x,y)f(y)dy$$
with $\Re(z)=a_{\beta}$. Observe that $d\geq k-C$ for some constant $C>0$. It
is convenient to divide the summation $\sum W_{k,l,d}$ over $d\geq k-C$ into two parts $|k-d|\leq N_1$ and $d>k+N_1$ with $N_1$ sufficiently large. More precisely, we may assume $N_1$ is so large that
$|y-\alpha_{i}x|\approx 2^{-k}$ in the support of $W_{k,l,d}$ for $2\leq i \leq m$. It is also clear that $W_{k,l,d}W_{k',l',d'}^{\ast}=0$ for $|l-l'|\geq 2$. Assume $|l-l'|\leq 1$. Then
\begin{eqnarray*}
\|W_{k,l,d}W_{k',l',d'}^{\ast}\|
\leq
C2^{-(d+d')}\left(2^{-dm_1}2^{-k(n-2-m_1-\beta)}\right)^{a_\beta}
\left(2^{-d'm_1}2^{-k(n-2-m_1-\beta)}\right)^{a_\beta}
\end{eqnarray*}
and $\|W_{k,l,d}W_{k',l',d'}^{\ast}\|$ is also less than or equal to
\begin{eqnarray*}
C_z\left(|\lambda|2^{-(d\wedge d')m_1}2^{-k(n-2-m_1)}\right)^{-1}\left(2^{-dm_1}2^{-k(n-2-m_1-\beta)}\right)^{a_\beta}
\left(2^{-d'm_1}2^{-k(n-2-m_1-\beta)}\right)^{a_\beta}.
\end{eqnarray*}
with $C_z=C(1+|\Im(z)|)^4$. Taking a convex combination, we obtain
\begin{eqnarray*}
\|W_{k,l,d}W_{k',l',d'}^{\ast}\|
&\leq& C_z|\lambda|^{-\frac{1}{\beta+1}}
2^{-|d-d'|\delta}\\
&\leq & C(1+|\Im(z)|)^4|\lambda|^{-\frac{1}{\beta+1}}2^{-|d-d'|\delta}
\end{eqnarray*}
with $\delta=\beta/(\beta+1)+m_1a_\beta>0$.
A similar argument gives us that $\|W_{k,l,d}^{\ast}W_{k',l',d'}\|$ is
also bounded by a constant multiple of
$|\lambda|^{-\frac{1}{\beta+1}}2^{-|d-d'|\delta}$ with
$\delta=\beta/(\beta+1)+m_1a_\beta>0$. By the
almost orthogonality principle, it follows that
$$\left\|\sum W_{k,l,d}\right\|\leq C(1+|\Im(z)|)^2
|\lambda|^{-\frac{1}{2(\beta+1)}}$$
where the summation is taken over all $k,l,d\geq 0$ satisfying
$|k-l|\leq N_0$ and $d\geq k+N_1$.

For $1\leq \omega \leq m$, we use $W_{k,l,d_1,\cdots,d_{\omega}}$ to
denote the operator obtained from $W_{k,l,d_1,\cdots,d_{\omega-1}}$ by
multiplying a factor $\Phi_{d_{\omega}}(y-\alpha_{\omega}x)$ in the
cut-off of $W_{k,l,d_1,\cdots,d_{\omega-1}}$.
Generally, we may apply the same argument as above to obtain
\begin{equation*}
\left\|\sum W_{k,l,d_1,\cdots,d_{\omega}}\right\|
\leq
C(1+|\Im(z)|)^2|\lambda|^{-\frac{1}{2(\beta+1)}},
\qquad 1\leq \omega \leq m
\end{equation*}
where the summation is taken over all nonnegative integers $k,l,d_1,
\cdots,d_{\omega}$ satisfying $|k-l|\leq N_0$, $|d_i-k|\leq N_i$ with
$1\leq i \leq \omega-1$ and $d_{\omega}\geq k+N_{\omega}$, and
$N_1,\cdots,N_{\omega}$ are sufficiently large
integers depending on the parameters $\alpha_i$ appearing in the
factorization of $S''_{xy}$. Thus it remains to show that
\begin{equation*}
\left\|\sum W_{k,l,d_1,\cdots,d_{m}}\right\|
\leq
C(1+|\Im(z)|)^2|\lambda|^{-\frac{1}{2(\beta+1)}}
\end{equation*}
for $|k-l|\leq N_0$ and $|k-d_i|\leq N_i$ with $1\leq i \leq m$.
Observe that
$W_{k,l,d_1,\cdots,d_{m}}W_{k',l',d_1',\cdots,d_{m}'}^{\ast}=0$ if $|l-l'|\geq 2$ and $W_{k,l,d_1,\cdots,d_{m}}^{\ast}W_{k',l',d_1',\cdots,d_{m}'}=0$ if $|k-k'|\geq 2$. By the almost orthogonality lemma, it is enough to show
that $\|W_{k,l,d_1,\cdots,d_{m}}\|$ is uniformly bounded above
by $C_z|\lambda|^{-\frac{1}{2(\beta+1)}}$. Indeed, we have
\begin{eqnarray*}
\left\|W_{k,l,d_1,\cdots,d_{m}}\right\|
&\leq& C_z\min\left\{\left(|\lambda|2^{-k(n-2)}\right)^{-1/2}
2^{-k(n-2-\beta)a_{\beta}},
 2^{-k}2^{-k(n-2-\beta)a_{\beta}}\right\}\\
&\leq & C(1+\Im|z|)^2\lambda|^{-\frac{1}{2(\beta+1)}}.
\end{eqnarray*}
The proof of the theorem is complete.
\end{proof}

As pointed out at the beginning of this section, Theorem
\ref{OIO with damped factor on L2} should be adapted in some cases for our purpose. Now we turn our attention to the special case when $S''_{xy}$ is of the form $cy^{\beta}(y-\alpha x)^{n-2-\beta}$. In the next section, we will see that Theorem \ref{OIO with damped factor on L2} is not suitable for
interpolation. Thus the damping factor shall be replaced by another one.

\begin{theorem}\label{section 5. theorem new damped OIO}
Assume $S$ is a real-valued homogeneous polynomial in two variables with
degree $n\geq 3$. If the Hessian $S''_{xy}$ equals
$cy^{\beta}(y-\alpha x)^{n-2-\beta}$
with $0<\beta <(n-2)$ and nonzero numbers $c$ and $\alpha$, then the
conclusion of Theorem {\rm\ref{OIO with damped factor on L2}} is also valid
with the damping factor in (\ref{section 5. standard damping factor})
replaced by
\begin{equation}\label{section 5. new damped factor}
D(x,y)=x(y-\alpha x)^{n-3-\beta}.
\end{equation}
\end{theorem}

\begin{proof}
Our argument is the same as that of above theorem. Let
$W_{\lambda}$ be defined as in (\ref{OIO with damped factor and cut-off})
with the new damping factor $D(x,y)$. We can decompose $W_{\lambda}$ as $\sum_{k,l\geq 0}W_{k,l}$ with $W_{k,l}$ given by (\ref{section 5. def of W(k,l)}). By dilation, we may assume $\alpha=1$.

 {\bf Estimates in the range $k\geq l+10$.}\\
Assume $|l-l'|\leq 1$. The size estimates of $W_{k,l}W_{k',l'}^{\ast}$
are given by
\begin{eqnarray*}
\|W_{k,l}W_{k',l'}^{\ast}\|
&\leq & C2^{-(k+l)/2}2^{-(k'+l')/2}
\left(2^{-k}2^{-l(n-3-\beta)}\right)^{a_{\beta}}
\left(2^{-k'}2^{-l'(n-3-\beta)}\right)^{a_{\beta}}
\end{eqnarray*}
and the oscillatory estimates assert that $\|W_{k,l}W_{k',l'}^{\ast}\|$
is bounded by
\begin{equation*}
 C(1+|\Im(z)|)^{4}
\left(|\lambda|2^{-l(n-2)}\right)^{-1}
\left(2^{-k}2^{-l(n-3-\beta)}\right)^{a_{\beta}}
\left(2^{-k'}2^{-l'(n-3-\beta)}\right)^{a_{\beta}}.
\end{equation*}
By a convex combination of above estimates, we obtain
$$\|W_{k,l}W_{k',l'}^{\ast}\|\leq C_z|\lambda|^{-\frac{1}{\beta+1}}
2^{-|k-k'|\delta}$$
with $\delta=\beta/[2(\beta+1)]+a_{\beta}>0$. In fact, we have
$$\delta=\frac1{2(\beta+1)}\left(\beta+1-\frac{\beta}{n-2-\beta}\right)>0.$$
Similarly, it is also true that
\begin{equation*}
\|W_{k,l}^{\ast}W_{k',l'}\|\leq C_z|\lambda|^{-\frac{1}{\beta+1}}
2^{-|l-l'|\delta}
\end{equation*}
with $|k-k'|\leq 1$ and $\delta=a_{\beta}+(n+\beta-2)/[2(\beta+1)]>0$.

 {\bf Estimates in the range $l\geq k+10$.}\\
 For $l\geq k+10$ and $l'\geq k'+10$ with $|l-l'|\leq 1$, we also have
 \begin{eqnarray*}
 \|W_{k,l}W_{k',l'}^{\ast}\|
 &\leq& C2^{-(k+l)/2}2^{-(k'+l')/2} 2^{-(k+k')(n-2-\beta)a_{\beta}}\\
 \|W_{k,l}W_{k',l'}^{\ast}\|
 &\leq& C_z\left(|\lambda|2^{-l\beta}2^{-(k\wedge k')(n-2-\beta)}\right)^{-1}
 2^{-(k+k')(n-2-\beta)a_{\beta}}.
 \end{eqnarray*}
Then it follows that $\|W_{k,l}W_{k',l'}^{\ast}\|$ is not greater than
a constant multiple of
$|\lambda|^{-\frac{1}{\beta+1}}2^{-|k-k'|\delta}$ with
$\delta=(n-\beta-2)/[2(\beta+1)]$. A similar argument shows that
$\|W_{k,l}^{\ast}W_{k',l'}\|$ is bounded by
$C_z|\lambda|^{-\frac{1}{\beta+1}}2^{-|l-l'|\delta}$ with
$\delta=\beta/[2(\beta+1)]>0$ for $|k-k'|\leq 1$.

  {\bf Estimates in the range $|k-l|\leq 10$.}\\
A further decomposition of $T_{k,l}$ is necessary to separate $(x,y)$
from the line $x=y$ on which $S''_{xy}$ vanishes. Let $T_{k,l,m}$ be
defined as $T_{k,l}$ with the cut-off of $T_{k,l}$ multiplied by the factor
$\Phi_m(x-y)$. It is clear that $m$ is not less than $k-12$. The
treatment in this case is more direct and does not need the almost
orthogonality principle. Actually, we shall see that the summation
$\sum W_{k,l,m}$ is absolute convergent. The size of the support implies
\begin{equation*}
\|W_{k,l,m}\|
\leq C 2^{-m}\left(2^{-k}2^{-m(n-\beta-3)}\right)^{a_{\beta}}.
\end{equation*}
The oscillatory estimate gives
\begin{equation*}
\|W_{k,l,m}\|
\leq C (1+|\Im(z)|)^2
\left(|\lambda|2^{-k\beta}2^{-m(n-\beta-2)}\right)^{-1/2}
\left(2^{-k}2^{-m(n-\beta-3)}\right)^{a_{\beta}}.
\end{equation*}
When $|\lambda|\geq 2^{k\beta}2^{m(n-\beta)}$, we take the oscillatory
estimate as the bound of $\|W_{k,l,m}\|$. Observe that $m\geq k-12$ in
the present situation. For those $k$ satisfying
$2^{k\beta}2^{(k-12)(n-\beta)}\leq |\lambda|$, we use $m_{k,\lambda}$
to denote the largest integer $m$ satisfying
$2^{k\beta}2^{m(n-\beta)}\leq |\lambda|$. By this definition, it is easy
to see that $|\lambda|\approx 2^{k\beta}2^{m_{k,\lambda}(n-\beta)}$. Now
we observe that
\begin{eqnarray*}
\sum_{m \leq m_{k,\lambda}}\|W_{k,l,m}\|
&\leq& C_z \sum_{ m \leq m_{k,\lambda}}
\left(|\lambda|2^{-k\beta}2^{-m(n-\beta-2)}\right)^{-1/2}
\left(2^{-k}2^{-m(n-\beta-3)}\right)^{a_{\beta}}\\
&\leq& C_z
\left(|\lambda|2^{-k\beta}2^{-m_{k,\lambda}(n-\beta-2)}\right)^{-1/2}
\left(2^{-k}2^{-m_{k,\lambda}(n-\beta-3)}\right)^{a_{\beta}}\\
&\leq & C_z
|\lambda|^{-1/2}2^{k(\beta/2-a_{\beta})}
\left(|\lambda|2^{-k\beta}\right)^{\delta_{\beta}}
\end{eqnarray*}
with
$\delta_{\beta}=\frac{1}{n-\beta}\left(\frac{n-2-\beta}{2}
-(n-\beta-3)a_{\beta}\right).$
Recall that $2^{k\beta}2^{(k-12)(n-\beta)}\leq |\lambda|$. Hence
\begin{eqnarray*}
\sum_{2^{kn}\leq C|\lambda|}
\sum_{m=k-13}^{m_{k,\lambda}}\|W_{k,l,m}\|
&\leq & C_z\sum_{2^{kn}\leq C|\lambda|}
|\lambda|^{-1/2}2^{k(\beta/2-a_{\beta})}
\left(|\lambda|2^{-k\beta}\right)^{\delta_{\beta}}\\
&\leq &C(1+|\Im(z)|)^2|\lambda|^{-\frac{1}{2(\beta+1)}},
\end{eqnarray*}
where the exponent of $2^k$ equals
\begin{equation}\label{section 5. reason for crucial assumption}
\frac{n}{n-\beta}\frac{\beta}{2(\beta+1)}-\frac{n}{n-\beta}
\frac{n-2(\beta+1)}{2(\beta+1)}\frac{1}{n-\beta-2}>0
\end{equation}
in the summation over $2^{kn}\leq C|\lambda|$. We must stress that the
assumption $\beta>0$ is crucial in the inequality
(\ref{section 5. reason for crucial assumption}). Similarly, we also
have
$$\sum_{2^{kn}\leq C|\lambda|}
\sum_{m=m_{k,\lambda}}^{\infty}\|W_{k,l,m}\|
\leq C|\lambda|^{-\frac{1}{2(\beta+1)}}.$$

For those large $k$ satisfying
$2^{k\beta}2^{(k-12)(n-\beta)}>|\lambda|$, we use the size estimate
as the upper bound of $\|W_{k,l,m}\|$. Then the desired estimate
follows directly. Actually, the size estimate gives
\begin{eqnarray*}
\sum_{2^{kn}\geq C|\lambda|}\sum_{m\geq k-12}\|W_{k,l,m}\|
&\leq& C\sum_{2^{kn}\geq C|\lambda|}\sum_{m\geq k-13}
2^{-m}\left(2^{-k}2^{-m(n-\beta-3)}\right)^{a_{\beta}}\\
&\leq & C\sum_{2^{kn}\geq C|\lambda|}2^{-kn/[2(\beta+1)]}\\
&\leq &C |\lambda|^{-\frac{1}{2(\beta+1)}}.
\end{eqnarray*}
The proof is therefore complete.
\end{proof}

\noindent{\bf Remark.} When $\beta=0$, the inequality
(\ref{section 5. reason for crucial assumption}) is not true and
its left side equals $-\frac{1}{2}$. Hence the argument above fails.
We may ask whether Theorem \ref{section 5. theorem new damped OIO} is
still true in the case $\beta=0$. The answer is negative.
Actually, we would obtain that the following operator
\begin{equation*}
Wf(x)=\int_{-\infty}^{\infty}e^{i(x-y)^n}|x(x-y)^{n-3}|^{1/2}f(y)dy
\end{equation*}
is bounded from $L^2$ to itself with $n\geq 3$ if Theorem
\ref{section 5. theorem new damped OIO} were still true when $\beta=0$.
Let $f$ be the characteristic function of the interval
$(N,N+\sqrt[n]{\pi/16})$ for $N\geq 1$. Then we see that
$|Wf(x)|\geq CN^{1/2}$ for
$N+\sqrt[n]{\pi/8}\leq x \leq N+\sqrt[n]{\pi/4}$. This observation
implies that $W$ is unbounded on $L^2$. The following theorem is a
substitute of the above theorem.

\begin{theorem}\label{section 5. translation invariant OIO}
Let $U$ be the operator as in Theorem {\rm\ref{H1 boundedness of
trans-inv OIO}} with $\Re(z)=(n-2)/2$. Then $U$ extends as a bounded operator on $L^{2}$ with the operator norm not greater than a constant multiple of $(1+|\Im(z)|)^{2}$.
\end{theorem}
\noindent{\bf Remark.} There is an equivalent formulation of the
theorem. Assume $S(x,y)=c(x-\alpha y)^n$ for nonzero real numbers $c$
and $\alpha$. Let $W_{\lambda}$ be defined as in
(\ref{OIO with damped factor and cut-off}) with $\Re(z)=(n-2)/2$ and the
damping factor
\begin{equation}\label{section 5. damping factor for trans-inv}
D(x,y)=(|\lambda|^{-1/n}+|x-\alpha y|).
\end{equation}
Then the estimate for $W_{\lambda}$ in
(\ref{decay estimate of OIO with Damped factor}) is still true.\\

\begin{proof}
Assume $c=\alpha=1$ without loss of generality. By the remark, it suffices to show that $W_{\lambda}$ satisfies the estimate (\ref{decay estimate of OIO with Damped factor}). Let $W_{k,l}^{\sigma_1,\sigma_2}$ be defined as in (\ref{section 5. def of W(k,l)}). In the case $k\geq l+10$, let $U_l=\sum(W_{k,l}^{++}+W_{k,l}^{-+})$ with the summation being taken over
$k\geq l+10$. It is easy to see that $U_l$ is supported in the rectangle $|x|\leq C2^{-l}$ and $y\approx 2^{-l}$. By the Schur lemma, we obtain
$$\sum\|U_l\|
\leq C\sum 2^{-l}|\lambda|^{-(n-2)/(2n)}
\leq C|\lambda|^{-1/2}$$
where the summation is taken over $2^{ln}\geq |\lambda|$. For those $l$
satisfying $2^{ln}<|\lambda|^{1/n}$, $U_lU_{l'}^{\ast}=0$ for
$|l-l'|\leq 1$ and it follows from Lemma
\ref{section 5. orthogonality structure} that $\|U_l^{\ast}U_{l'}\|$ is
bounded by a constant multiple of
$(1+|\Im(z)|)^4|\lambda|^{-1}2^{-|l-l'|(n-2)/2}$. The
estimate in the case $l\geq k+10$ can be treated similarly. Actually,
set $V_k=\sum_{l\geq k+10}(W_{k,l}^{++}+W_{k,l}^{+-})$ for each $k$.
Then we can divide the summation $\sum V_k$ into two parts, i.e.
$2^{kn}\geq |\lambda|$ and $2^{kn}\leq |\lambda|$. The desired estimate
follows by using the Schur lemma to the first part and the almost
orthogonality to the second one.

 Now consider the case $|k-l|\leq 10$. Define $W_{k,l,m}$ as $W_{k,l}$
with the cut-off of $W_{k,l}$ multiplied by the factor $\Phi_m(x-y)$.
 Then $W_{k,l,m}=0$ unless $m\geq k-13$. For those $m$ satisfying
 $2^{mn}>|\lambda|$, let $U_m=\sum_{|k-l|\leq 10}W_{k,l,m}$ for each
 $m$. Then by Schur's lemma we obtain
 $\sum \|U_m\|\leq C|\lambda|^{-1/2}$ with the summation taken over
 all $m$ satisfying $2^{mn}>|\lambda|$. For $m$ satisfying
 $2^{mn}\leq |\lambda|$, we shall use the almost orthogonality. Indeed,
 we have
 \begin{eqnarray*}
 \|W_{k,l,m}W_{k',l',m'}^{\ast}\|
 &\leq &C(1+|\Im(z)|)^4|\lambda|^{-1}2^{-|m-m'|(n-2)/2},
 \quad |l-l'|\leq 1\\
  \|W_{k,l,m}^{\ast}W_{k',l',m'}\|
  &\leq &C(1+|\Im(z)|)^{4}|\lambda|^{-1}2^{-|m-m'|(n-2)/2},
 \quad |k-k'|\leq 1.
 \end{eqnarray*}
Thus we have almost orthogonality which implies the desire estimate.
\end{proof}

The following theorem implies the $L^2$ boundedness of $U$ in Theorem
\ref{H1 boundedness of trans-inv OIO}.
\begin{theorem}
Assume $n\geq 2$ is an integer. Let $W_{\lambda}$ be given by
\begin{equation}
W_\lambda f(x)=\int_{-\infty}^{\infty}e^{i\lambda (x-y)^n}
\left(|\lambda|^{-1/n}+|x-y|\right)^z \varphi(x,y)f(y)dy
\end{equation}
with $\Re(z)=-1$ and $\varphi \in C_0^\infty$. Then the operator norm
$\|W_\lambda\|$ on $L^2$ is bounded by a constant multiple of $(1+|\Im(z)|)^2$ with the
constant independent of $\lambda$.
\end{theorem}
\begin{proof}
Define $W_{k,l}^{\sigma_1,\sigma_2}$ as in the proof of Theorem \ref{OIO with damped factor on L2}. Then we shall divide the proof into three cases $k\geq l+10$, $l\geq k+10$ and $|k-l|\leq 10$.
In the case $k\geq l+10$, define $U_l=\sum(W_{k,l}^{++}+W_{k,l}^{-+})$ for each $l$. For those $l$ satisfying $2^{ln}\geq |\lambda|$, we apply Schur's lemma to obtain
$$\left\|\sum U_l\right\|\leq C\sum_{2^{ln}
\geq |\lambda|}2^{-l}|\lambda|^{1/n}\leq C.$$
For $l$ and $l'$ satisfying $2^{ln}\leq |\lambda|$ and
$2^{l'n}\leq |\lambda|$, we invoke Lemma \ref{section 5. orthogonality structure} to obtain
$$\|U_l^{\ast}U_{l'}\|
\leq C_z (|\lambda|2^{-(l\wedge l')(n-2)})^{-1}2^{l+l'}
\leq C_z 2^{-|l-l'|(n-1)}$$
with $C_z$ independent of $\lambda$. Observe $U_lU_{l'}^{\ast}=0$ for
$|l-l'|\geq 2$. We have the almost orthogonality and then obtain
$\|\sum U_l\|\leq C$. The treatment of the case $l\geq k+10$ is similar.
For $|k-l|\leq 10$, we shall introduce $W_{k,l,m}$ as above. For those
$m$ satisfying $2^{mn}\geq |\lambda|$, the desired estimate follows by
the Schur's lemma. While for $2^{mn}\leq |\lambda|$, it is true that
$$ \|W_{k,l,m}W_{k',l',m'}^{\ast}\|
 \leq C(1+|\Im(z)|)^42^{-|m-m'|(n-1)},\quad |l-l'|\leq 1.$$
The same upper bound is also valid for $\|W_{k,l,m}^{\ast}W_{k',l',m'}\|$
for $|k-k'|\leq 1$. By the almost orthogonality principle, we obtain
$\left\|\sum W_{k,l,m}\right\|\leq C.$ The proof is complete.
\end{proof}

\section{Proof of Theorem \ref{main theorem 1}}\label{proof of the main result}

\hspace{.5cm}In this section, we shall apply previous results to prove Theorem \ref{main theorem 1}.
Let $p_{0}$ and $p_{1}$ be equal to $n/(n-k_{min})$ and $n/(n-k_{max})$, respectively. The argument will be divided into two cases: (i) $k_{min}<n/2$ and (ii) $k_{min}\geq n/2$. In both cases, there may occur $k_{min}=k_{max}$, i.e., $S$ is a monomial modulo pure-$x$ and pure-$y$ terms. We shall see that the following argument in the case (ii) is also applicable to the case when $S$ is a monomial.

\begin{proof}
By duality, it is enough to show that $T$ is bounded on $L^{p_0}$. If this were done, we would obtain that its adjoint operator $T^\ast$ is bounded from $L^{n/k_{max}}$ to itself. Thus $T$ has a bounded extension from $L^{p_1}$ to itself.

We first treat the case (ii) $k_{min}\geq n/2$. Since $T^\ast$ has the oscillating kernel $e^{-iS(y,x)}$, we see that the assumption $k_{min}\geq n/2$ imposed on $T$ is equivalent to the condition $k_{max}\leq n/2$ on $T^\ast$. Hence it suffices to prove that $T$ is bounded from $L^{n/(n-k_{max})}$ to itself under the assumption $k_{max}\leq n/2$. For brevity, we only prove that
$T$ maps $L^{n/(n-k_{max})}(0,\infty)$ to itself boundely since other parts of integration can be treated similarly.

Consider
\begin{equation*}
H(f)(x)=\int_{0}^{\infty}e^{iS\l(x^{k_{max}/(n-k_{max})},y\r)}f(y)dy
\end{equation*}
and denote by $H^\ast$ its adjoint operator. By the van der Corput lemma, we have
\begin{eqnarray*}
\int_{0}^{\infty}\left|H^\ast(f)(x)\right|^2dx
 &\leq &
C|a_{n-k_{max}}|^{-1/k_{max}}\int_{0}^{\infty}\int_{0}^{\infty}
\left|y_{1}^{k_{max}}-y_{2}^{k_{max}}\right|^{-\frac1{k_{max}}}|f(y_{1})f(y_{2})|dy_{1}dy_{2}\\
& \leq & C|a_{n-k_{max}}|^{-1/k_{max}}\|f\|_{2}^{2},
\end{eqnarray*}
where we impose the assumption on $k_{max}>1$. Otherwise if $k_{max}=1$, then $H$ reduce to the Fourier transform. Hence $H$ is bounded on $L^2(0,\infty)$. By Lemma \ref{interpolation lemma}, we have
\begin{equation*}
\int_0^\infty |H(f)(x)|^{p}x^{p-2}dx
\leq C\|f\|_{p}^{p}
\end{equation*}
for $1< p \le 2$. Set $p=n/(n-k_{max})$. By a change of variables, we see that
$T$ is bounded on $L^{p}(0,\infty)$. When $k_{min}=k_{max}$, the range of $p$ described in Theorem \ref{main theorem 1} is just a single value $p=n/(n-k_{min})$. If $k_{min}\geq n/2$, then we have proved that $T$ is bounded on $L^p$. If $k_{min}<n/2$, we can apply a duality argument to obtain the desired statement.

Now we turn to the case (i) $1\le k_{min}<\frac{n}{2}$. Without loss of generality, we assume also $k_{min}<k_{max}$. The argument is somewhat different depending
on whether the Hessian $S''_{xy}$ is of the form $cy^{\beta}(y-\alpha x)^{n-\beta-2}$ with nonzero real numbers $c$ and $\alpha$. We shall divide the argument into three cases. Recall that $S_{xy}''$ can be written as
$$ S''_{xy}(x,y)=Cx^\gamma
y^\beta\prod\limits^m_{j=1}(y-\alpha_jx)^{m_j}\prod\limits^s_{j=1}Q_j(x,y).$$
Since $k_{min}<n/2$, it is easy to see $0\leq \beta <(n-2)/2$.

Case I: $S_{xy}''$ is not of the form $cy^{\beta}(y-\alpha x)^{n-\beta-2}$.
Then let $D$ be the damping factor given by
\begin{equation}\label{sec6. damped factor one}
D(x,y)=x^\gamma\prod\limits^m_{j=1}(y-\alpha_jx)^{m_j}\prod\limits^s_{j=1}Q_j(x,y).
\end{equation}
 Consider the following analytic family of operators
\begin{equation}\label{sec6. analytic family operators}
T_{z}f(x)=\int_{-\infty}^{\infty}e^{iS(x,y)}|D(x,y)|^{z}f(y)dy
\end{equation}
for $z$ in the strip $\l\{x+iy:-\frac{1}{n-2-\beta}\le x \le
a_{\beta}\r\}$. It is clear that $T_0=T$. Let $K_{z}(x,y)=|D(x,y)|^{z}.$
Choose $\alpha_{j} \neq 0$ for $m+1 \le j \le m+2s$ such that $\alpha_{j}\neq \alpha_{k}$ for
$1\le j\ne k \le m+2s$. For $\Re(z)=-1/(n-2-\beta)$, the kernels $K_{z}$
fall under the scope of the class discussed in Theorem
\ref{Fraction integrals on Hardy spaces} and Theorem
\ref{osc fractional operators with polynomial phases}. Indeed,
for $\Re(z)=-1/(n-2-\beta)$, a direct computation shows that
$$ |K_{z}(x,y)|\le C|x|^{-\theta_{0}}\prod_{1\le j \le
 m+2s}|x-c_{j}y|^{-\theta_j},$$
$$|\partial_{y}K_{z}(x,y)|\le C|x|^{-\theta_{0}}\sum_{1\le k \le
m+2s}|x-c_{k}y|^{-\theta_{k}-1}\prod_{j\ne k}|x-c_{j}y|^{-\theta_j},
$$
  where
$$\theta_{0}=\frac{\gamma}{n-2-\beta},\quad\theta_{j}=\frac{m_j}{n-2-\beta},\;\, 1\le j \le m\;\;
\textrm{and}\;\;
\theta_{j}=\frac{1}{n-2-\beta},\;\,m+1\le j \le m+2s$$ with
$c_{j}=\alpha_{j}^{-1}\neq 0$ for $1\le j \le m+2s$.
On the one hand, it follows from
Theorem \ref{osc fractional operators with polynomial phases}
that $T_{z}$ is bounded from $H_{E}^{1}$ to $L^{1}$ with the bound
less than a constant multiple of $(1+|\Im(z)|)^{2}$.
On the other hand, it follows from
Theorem \ref{OIO with damped factor on L2} that $T_{z}$ is bounded on
$L^{2}$ when $\Re(z)=a_{\beta}$ with the bound not greater than
$C(1+|\Im(z)|)^2$. Combing these results, we obtain
$$\|T_{z}^{*}f\|_{BMO_E} \leq A_{z}\|f\|_{L^{\infty}}$$
for $ \Re(z)=-1/(n-2-\beta)$, and $$  \|T_{z}^{*}f\|_{L^{2}} \leq
A_{z}\|f\|_{L^{2}}$$ for $\Re(z)=a_{\beta}$, where $BMO_{E}$ is given by Definition \ref{definition of variants of BMO} with $P(x,y)=-S(y,x)$. It follows immediately that
\begin{equation}\label{sec6. L infty estimate}
\l\|\big(T_z^{\ast}f\big)^{\sharp}_E\r\|_{L^{\infty}} \leq
A_{z}\|f\|_{L^{\infty}}
\end{equation}
for $ \Re(z)=-1/(n-2-\beta)$, and
\begin{equation}\label{sec6. L2 estimate}
\l\|\big(T_{z}^{\ast}f\big)^{\sharp}_E\r\|_{L^{2}}
\leq A_{z}\|f\|_{L^{2}}
\end{equation} for $ \Re(z)=a_{\beta}$.

To apply the complex interpolation, we need a linear operator to
approximate $(T_z^{\ast}f)^{\sharp}_E$. Let $\rho(x,y)$ be a measurable
function with $|\rho(x,y)|\leq 1$. We use $Q(x)$ to denote a measurable
mapping from $x$ to a cube $Q(x)$ containing $x$. Then we can define
a family $A_z$ of operators by
\begin{equation*}
A_zf(x)=\frac{1}{|Q(x)|}\int_{Q(x)}\left(T_z^{\ast}f(y)-
(T_z^{\ast}f)_{Q(x)}^E\right)
\rho(x,y)dy
\end{equation*}
with $(T_z^{\ast}f)_{Q(x)}^E$ defined as in the introduction. Then
$(T_z^{\ast}f)^{\sharp}_E$ can be obtained by taking the supremum
$\sup|A_zf(x)|$ over all $\rho(x,y)$ and $Q(x)$ described above.
By interpolation and Lemma \ref{interpolation of sharp function}, we
get that
\begin{equation}
\l\|T_{z}^{*}(f)\r\|_{L^p}\le
\l\|\big(T_{z}^{*}\big)^{\sharp}(f)\r\|_{L^{p}}\leq
C_{z}\|f\|_{L^{p}},
\end{equation}
where $p$ and $z$ satisfy
$$\Re(z)=-\frac1{n-2-\beta}(1-\theta)+\alpha_{\beta}\theta$$
and $1/p=\theta/2$ for $0\le \theta \le 1$. If $\Re(z)=0$, we obtain
$\theta=2(\beta+1)/n$ and $p=n/(\beta+1)$. By duality, $T_{0}$ is
bounded on $n/(n-(\beta+1))$. Since $k_{min}=\beta+1$, the desired
result follows.

Case II: $S_{xy}''=cy^{\beta}(y-\alpha x)^{n-\beta-2}$ with
$0<\beta<(n-2)/2$. Define the analytic family $T_z$ as in
(\ref{sec6. analytic family operators}) with the damping factor in
(\ref{sec6. damped factor one}) replaced by
$D(x,y)=x(y-\alpha x)^{n-\beta-3}.$
Then the above argument is also applicable here. Actually, when the real
part of $z$ is equal to $a_{\beta}$, $T_z$ is bounded
from $L^2$ to itself with the norm bounded by a constant multiple of
$(1+|\Im(z)|)^2$. For those $z$ with $\Re(z)=-1/(n-2-\beta)$, we may invoke Theorem
\ref{osc fractional operators with polynomial phases} to obtain that $T_z$ is bounded from $H_E^1$ to $L^1$. The interpolation argument is the same as above. Thus we show that $T$ is bounded on $L^{p_0}$.

Case III: $S_{xy}''=c(y-\alpha x)^{n-2}$. By dilation, we may assume
$\alpha=1$. The damping factor shall be replaced by
$D(x,y)=\left(1+|x-y|\right)^{n-2}.$
By combing Theorem \ref{H1 boundedness of trans-inv OIO} and
Theorem \ref{section 5. translation invariant OIO}, the desired result
follows by a similar argument.

The proof is therefore complete.
 \end{proof}

\section{Applications}

 \hspace{.5cm}In this section, we shall give two applications of Theorem \ref{main theorem 1}. We begin with the proof of Theorem \ref{Lp result for OIO with nonhomogeneous phases}.

\begin{proof}
By Theorem \ref{main theorem 1} and Lemma \ref{interpolation lemma},
we have
\begin{equation*}
\int_{0}^{\infty}\left|\int_{0}^{\infty}\exp(iS(x,y))f(y)dy\right|^{p}
x^{(p-q_{0})/(q_0-1)}dx \leq C \int_{0}^{\infty}|f(x)|^{p}dx
\end{equation*}
 for $n/(n-k_{min})\leq q_0 \leq n/(n-k_{max})$ and $1<p\leq q_0$.

  Set
 $$q_{0}=p_{0}=\frac{n}{n-k_{min}}
 \quad\;\textrm{and}\quad\; p=\widetilde{p_{0}}=\frac{(n-k_{min})m_1+k_{min}}{(n-k_{min})m_1}.$$
It follows that
\begin{equation}\label{the weighted boundedness of OIO with nonhomogeneous polynomials}
\int_{0}^{\infty}\left|\int_{0}^{\infty}\exp(iS(x^{m_1},y))f(y)dy\right|^{p}
dx \leq C \int_{0}^{\infty}|f(x)|^{p}dx
\end{equation}
with $p=\widetilde{p_0}$.
The above inequality also holds if we set
$$q_{0}=p_{1}=\frac{n}{n-k_{max}}
 \quad\;\textrm{and}\quad\;
p=\widetilde{p_{1}}=\frac{(n-k_{max})m_1+k_{max}}{(n-k_{max})m_1}.$$
 By interpolation, we obtain that the inequality (\ref{the weighted boundedness of OIO with nonhomogeneous polynomials})
 is still true for $\widetilde{p_0}\leq p \leq \widetilde{p_1}$.
 By a duality argument and the interpolation technique as above, we conclude that
 \begin{equation*}
\int_{0}^{\infty}\left|\int_{0}^{\infty}\exp\l(iS(x^{m_1},y^{m_2})\r)f(y)dy\right|^{p}
dx \leq C \int_{0}^{\infty}|f(x)|^{p}dx
\end{equation*}
 holds for $p$ in the range
 $$\frac{(n-k_{min})m_1+k_{min}m_2}{(n-k_{min})m_1}\leq p \leq \frac{(n-k_{max})m_1+k_{max}m_2}{(n-k_{max})m_1}.$$
 The other part of $T_{m_1,m_2}f$ can be treated similarly.
 Hence the $L^p$ boundedness of $T_{m_1,m_2}$ has been established.
 By dilation as in \cite{PS1997}(also \cite{yangchanwoo2}), we can show that
 this range is also sharp. The proof is therefore complete.
 \end{proof}

Pitt's inequality is a generalized Hausdorff-Young inequality with power weights; see Beckner \cite{beckner}. We shall give a simple proof by using previous results. Denote by $\widehat{f}$ the Fourier transform of $f$ in the Schwartz class,
\begin{equation*}
\widehat{f}(x)=\int_{\mathbb{R}^n}e^{-2\pi ix\cdot y}f(y)dy.
\end{equation*}
Then Pitt's inequality can be stated as follows. For $1\le p \le \infty$, let $p'$
be the conjugate exponent of $p$, i.e. $1/p+1/p'=1$.
\begin{theorem}\label{Theorem7.1}
{\rm(\cite{beckner})} Let $1<p\leq q<\infty$, $0\le \alpha<n/q$ and $0\le \beta <n/p'$. There exists a constant $C>0$ such that
\begin{equation}
\left(\int_{\mathbb{R}^n}\left||x|^{-\alpha}\widehat{f}(x)\right|^q dx\right)^{1/q}
\leq C\left(\int_{\mathbb{R}^n}\left||x|^\beta f(x)\right|^p dx\right)^{1/p}
\end{equation}
for all Schwartz functions $f$, where $p,~q,~\alpha$ and $\beta$ satisfy
$n/p+n/q+\beta-\alpha=n$.
\end{theorem}
\begin{proof}
We first prove the theorem in one dimension. For $n=1$,
by dividing the integration into two parts, it is enough to show that
\begin{equation}\label{OIO with monomial polynomial phase}
\left(\int_{0}^{\infty}\left|\int_0^\infty e^{-2\pi ixy}f(y)dy\right|^{q}x^{-q\alpha}dx\right)^{1/q}\leq
C \left(\int_{0}^{\infty}|f(x)|^{p}x^{p\beta}dx\right)^{1/p}.
\end{equation}
Consider
$$Tg(x)=\int_{0}^{\infty}\exp\l(-2\pi i x^ay^b\r)g(y)dy,$$
for $a,b\ge 1$. By the same argument as in $\S$\ref{proof of the main result} when $S$ is a monomial, we can show that $T$ is bounded from $L^{(a+b)/a}(0,\infty)$ to itself. By dilation and
interpolation, it is easy to see that $T$ has a bounded extension from $L^{p}$ to $L^{q}$ with $1\le p \le (a+b)/a$ and $1/p=1-b/(aq)$. Set $a=1/(1-q\alpha)$ and $b=1/(1-p'\beta)$. The inequality (\ref{OIO with monomial polynomial phase}) follows by a change of variables and letting $g(y)=y^{b-1}f\l(y^b\r).$

 For dimension $n \ge 2$, observe that
$|x|^{-\alpha}\le \prod_{k=1}^{n}|x_{k}|^{-\alpha/n}$ and
$|x|^{\beta}\ge \prod_{k=1}^{n}|x_{k}|^{\beta/n}$
for $x\in \mathbb{R}^n$ and $\alpha,\beta \ge 0$. Now we define $\mathscr{F}_{k}$ to be the Fourier transform relative to $x_{k}$ with other variables fixed. It is clear that $\mathscr{F}$ is just the composition of $\mathscr{F}_{k}$, i.e.,  $\mathscr{F}f(x)=\mathscr{F}_n\mathscr{F}_{n-1}\cdots\mathscr{F}_1f(x)$. By $n$ applications of the one dimensional inequality and Minkowski's inequality , we obtain Pitt's inequality in dimension $n$. Thus the proof is complete.
\end{proof}

\end{document}